\documentclass[12pt]{article}
\newcommand{\kmcomment}[1]{}
\usepackage[leqno]{amsmath}
\usepackage{amssymb,amsfonts}
\usepackage{mathtools}
\usepackage{newtxmath}

\newcommand{\kmqed}{\hfill\ensuremath{\blacksquare}}

\newcommand{\ovfrakg}{\overline{\frakg}}

\newcommand{\ds}{\ensuremath{\displaystyle }}
\newcommand{\pdel}{\partial}

\newcommand{\myHom}[1]{\textrm{H}_{#1}}

\newcommand{\mN}{\ensuremath{\mathbb{N}}} 
\newcommand{\mR}{\ensuremath{\mathbb{R}}} 
\newcommand{\mZ}{\ensuremath{\mathbb{Z}}} 

\newcommand{\frakg}{\mathfrak{g}}

\newcommand{\tbdl}[1]{\mathrm{T}(#1)}

\newcommand{\Sbt}[2]{[#1,#2]}                
\newcommand{\SbtS}[2]{[#1,#2]_{\text{\footnotesize S}}} 
                
\newcommand{\CSp}[1]{\text{C}_{#1}} 

\newcommand{\eps}[1]{\epsilon_{#1}}

\renewcommand{\dim}{\textrm{dim}}

\newcommand{\mywedge}{\bigtriangleup} 

\newcommand{\myCS}[1]{ \text{C}_{#1}} 

\renewcommand{\[}{$$} \renewcommand{\]}{$$}

\usepackage{ntheorem} 
{
\theoremheaderfont{\bfseries}
\theorembodyfont{\rmfamily}

\newtheorem{defn}{\textbf{Definition}}
\newtheorem{prop}{Proposition}[section]

\newtheorem{remark}{Remark}[section]

\newtheorem{definition}{Definition}[section]
\newtheorem{theorem}{Theorem}[section] 
\newtheorem{thm}{Theorem}[section]
\newtheorem*{thm-none}{Theorem}[section]

\newtheorem{lemma}[theorem]{Lemma}
\newtheorem{kmCor}[theorem]{Corollary}

}

\renewcommand{\[}{$$} \renewcommand{\]}{$$}

\newcommand{\wtedC}[2]{ C_{#1}^{[#2]}}      
\newcommand{\wtedCR}[2]{ \widetilde{C}_{#1}^{[#2]}}      
\newcommand{\wtedCRx}[2]{ \widehat{C}_{#1}^{[#2]}}      

\newcommand{\CYm}[2]{ \widehat{{\mathfrak Y}_{#2} }[#1]} %
\newcommand{\myL}[1]{\operatorname{L}_{#1}}
\newcommand{\myE}[2]{E^{#1}_{#2}}
\newcommand{\myWE}[2]{\widehat{E}^{#1}_{#2}}

\newcommand{\myF}[2]{\operatorname{F}^{#1}_{#2}}
\newcommand{\myTF}[2]{\widetilde{\text{F}}^{#1}_{#2}}
\newcommand{\Hori}[2]{\mathfrak{S}^{[#1]}_{#2}}

 \newcommand{\zb}[1]{z_{#1}}
 \newcommand{\ub}[1]{u_{#1}}
 \newcommand{\bU}[2]{{U}^{#1}_{#2}}
 \newcommand{\LC}[2]{\lambda^{#1}_{#2}}
 \newcommand{\Lc}[2]{{L}^{#1}_{#2}}
 \newcommand{\mikH}{H}
\newcommand{\mykappa}[1]{K_{#1}}

 \newcommand{\MC}[2]{{\mu}^{#1}_{#2}}
 \newcommand{\varE}{\mathcal{E}}
 \newcommand{\myUnital}[1]{\mathfrak{e}_{#1}}

 \newcommand{\SbtES}[2]{[#1,#2]_{res}}
 
\newcommand{\myeps}[1]{\tilde{\epsilon}_{#1}}

\newcommand{\HS}{\text{single star}}
\newcommand{\DS}{\text{twins}} 
\newcommand{\TS}{\text{triple stars}}

\usepackage{geometry}
\usepackage[all,warning]{onlyamsmath}
\numberwithin{equation}{section}

\usepackage{mathtools}

\usepackage{tikz}
\usetikzlibrary{positioning,calc, quotes, angles, decorations.pathmorphing}
\usetikzlibrary{backgrounds}

\makeatletter
\global\let\tikz@ensure@dollar@catcode=\relax
\makeatother

\usepackage{wrapfig}

\usepackage{amssymb,fge}

\title{Super homologies associated with low dimensional Lie algebras}

\author{Kentaro Mikami\thanks{Akita University }
 \and Tadayoshi Mizutani\thanks{Professor Emeritus, Saitama University }}

\allowdisplaybreaks
\AtBeginDocument{\parindent=0pt \setlength{\headheight}{16pt}} 
\begin{document}
\kmcomment{
\begin{titlepage} 
\vspace*{20mm}
\begin{center}
{\Large Super homologies associated with low dimensional Lie algebras}

\vspace{10mm}
{Kentaro Mikami ({Akita University}) 
 \quad Tadayoshi Mizutani ({Saitama University})}

\vspace{10mm}
{\today} 
\end{center}
\end{titlepage}
}
\maketitle
\tableofcontents

\section{Introduction}

A Poisson structure \(\pi\) on a manifold \(M\) is a 2-vector field
characterized as \(\SbtS{\pi}{\pi} = 0\), where \(\SbtS{\cdot}{\cdot}\) is
the Schouten bracket. The graded algebra 
\( \sum \Lambda^{p} \tbdl{M} \) with the Schouten bracket is a prototype of Lie
superalgebra (cf.\ \cite{Mik:Miz:super2}). The Poisson condition 
\(\SbtS{\pi}{\pi} = 0\) means that a 2-chain \( \pi \wedge \pi\) is a cycle. Thus, 
studying the second super homology group of Lie superalgebra of tangent
bundle of $M$  is an activity  of Poisson geometry
(cf.\ \cite{Mik:Miz:super3}).  Given a \(\mZ\)-graded Lie
superalgebra, the 0-graded subspace is a Lie algebra. In this note, 
using the DGA 
\( \sum \Lambda^{p} \tbdl{M} \) with the Schouten bracket as a model, we
start from an abstract Lie algebra, construct non-trivial Lie superalgebra 
by Schouten-like bracket. Then it may be natural to ask how the core Lie
algebra control the Lie superalgebra. One trial  here is to investigate the
Betti numbers of the super homology groups.  
For abelian Lie algebras, the boundary operator is trivial, and the Betti
number is equal to the dimension of chain space for a given weight.  
So, 
we study their super homology groups for low dimensional Lie algebras of
dimension smaller than 4.    

For non-abelian Lie algebras of dimension 2, we have 
\begin{center}
\(
 \begin{array}{c | *{5}{c}}
 w>0  &  \wtedC{w}{w} & \wtedC{w+1}{w} & \wtedC{w+2}{w} \\\hline
 \text{SpaceDim} & 1 & 2 & 1 \\\hline
 \text{KerDim} & 1 & 
 1 
  & 0 \\\hline
 \text{Betti} & 0 & 0 & 0 \\\hline
  \end{array}
 \)
\end{center}

We classify non-abelian Lie algebras \(\frakg\) of dimension 3 into four cases:

(1) \(\dim[\frakg,\frakg  ] =1 \) and \(\dim[\frakg,\frakg  ] \subset
Z(\frakg) \),   \qquad 
(2) \(\dim[\frakg,\frakg  ] =1 \) and \(\dim[\frakg,\frakg  ] \not\subset Z(\frakg) \),

(3) \(\dim[\frakg,\frakg  ] =2 \),  \qquad 
(4) \(\dim[\frakg,\frakg  ] =3 \).  

Even though the case (4) has two kinds of Lie algebras when the base field is
\(\mR\), we made the same treatment and got the same results.  

When \(w=0\) the
homology groups are the same with usual Lie algebra homology groups, so 
we show the
Betti numbers for super homology groups when \(w>0\). 

 
 \begin{center}
 \(
 \begin{array}{c|*{5}{c}}
 \text{weight}= w>0 & w-1 & w & w+1 & w+2 & w+3\\\hline
 \text{SpaceDim} & 
 \binom{w}{2} & 
 3 \binom{w}{2} + \binom{w+2}{2} & 
 3\binom{w}{2} + 3\binom{w+2}{2} & 
 \binom{w}{2} + 3\binom{w+2}{2} & 
 \binom{w+2}{2} 
 \\\hline\hline
 \text{(1)'s Betti} & 0 & w & 3w+1 & 3w+2  & w+1
 \\\hline
 \text{(2)'s Betti} & 0 & 1 & 2 & 1 & 0
 \\\hline
 \text{(3)'s Betti} & 0 & 0 & \kappa & 2 \kappa  & \kappa 
 \\\hline
 \text{(4)'s Betti} & 0 & 0 &  0  & 0   &  0 
 \end{array}
 \)
\end{center}
where \(\kappa\) in (3) is defined by 
\( \kappa = \begin{cases} 1 & \text{if\quad} \alpha = -1 \\ 0 & 
\text{if\quad} \alpha \ne  -1 \end{cases} \), 
and \(\alpha\) is a non-zero parameter of the Lie bracket relations
\(\Sbt{\zb{1}}{\zb{2}}= 0\), \(\Sbt{\zb{1}}{\zb{3}}= \zb{1}\), \(\Sbt{\zb{2}}
{\zb{3}}= \alpha \zb{2}\).

\kmcomment{
We get the full table of super homology groups for each Lie algebra of
\(\dim \leqq 3\), which are reported in corresponding subsection. For
instance, 3-dimensional Lie algebra \(\frakg\) with \( \dim [\frakg,\frakg]
=3\), then Betti numbers are all zero in Theorem \ref{thm:G3D3-T1}. 
}

\section{Preliminaries, Notations and Basic Facts} 
\begin{definition}[Lie superalgebra]
Suppose a real vector space 
$\frakg $ is graded by \(\ds \mZ\) as 
\(\ds \frakg = \sum_{j\in \mZ} \frakg_{j} \)
and has a bilinear operation \(\Sbt{\cdot}{\cdot}\)  satisfying 
\begin{align}
& \Sbt{ \frakg_{i}}{ \frakg_{j}} \subset \frakg_{i+j} \label{cond:1} \\
& \Sbt{X}{Y} + (-1) ^{ x y} 
 \Sbt{Y}{X} = 0 \quad \text{ where }  X\in \frakg_{x} \text{ and } 
 Y\in \frakg_{y}  \\
& 
(-1)^{x z} \Sbt{ \Sbt{X}{Y}}{Z}  
+(-1)^{y x} \Sbt{ \Sbt{Y}{Z}}{X}  
+(-1)^{z y} \Sbt{ \Sbt{Z}{X}}{Y}  = 0\;.  
\label{super:Jacobi}
\end{align}
Then we call \(\frakg\) a \textit{\(\mZ\)-graded} (or \textit{pre}) Lie superalgebra.    

\kmcomment{
A Lie  \textit{superalgebra} $\frakg $ is graded by \(\ds \mZ_{2}\) as \(\ds
\frakg = \frakg_{[0]} \oplus \frakg_{[1]} \) and the condition
\eqref{cond:1} is regarded as  \(\ds \Sbt{ \frakg_{[1]}}{ \frakg_{[1]}}
\subset \frakg_{[0]}\) in modulo 2 sense.   
}
\end{definition}

\begin{remark} 
Super Jacobi identity \eqref{super:Jacobi} above is equivalent to the one of
the following.  
\begin{align}
\Sbt{ \Sbt{X}{Y}}{Z} &=  
\Sbt{X}{ \Sbt{Y}{Z} } + (-1)^{y z} \Sbt{ \Sbt{X}{Z}}{Y} \\  
\Sbt{X}{\Sbt{Y}{Z}} &=  
\Sbt{\Sbt{X}{Y}}{Z}  + (-1)^{x y} \Sbt{Y}{ \Sbt{X}{Z}} 
\end{align} 
\end{remark}

As in a usual Lie algebra homology theory, $m$-th chain space is the
exterior product \(\ds \Lambda^{m}\frakg\) of \(\frakg\) and the
boundary operator essentially comes from the operation \(\ds X \wedge
Y \mapsto \Sbt{X}{Y}\), 
in the case of pre Lie superalgebras, 
"exterior algebra" is defined as the quotient of the tensor
algebra \(\ds \otimes^{m}  \frakg \) of \(\frakg \) by the two-sided
ideal generated by  \begin{align} & X \otimes Y + (-1)^{x y} Y \otimes X
\quad\text{where }\quad X \in \frakg_{x}, Y \in \frakg_{y} \;,
\end{align} and we denote the equivalence class of \(\ds X \otimes Y \)
by \(\ds X \mywedge Y\).  
\kmcomment{ \(\ds \bigtriangleup \triangle
\mywedge  \nabla \bigtriangledown \) } 
Since \( \ds X_{\text{odd}}
\mywedge Y_{\text{odd}}  = Y_{\text{odd}} \mywedge X_{\text{odd}} \) and
\( \ds X_{\text{even}} \mywedge Y_{\text{any}}  = - Y_{\text{any}}
\mywedge X_{\text{even}} \) hold, \( \mywedge ^{m} \frakg_{k} \) has a
symmetric property for odd $k$ and has a skew-symmetric property 
for even $k$ with respect to \(\mywedge\).

\medskip

Suppose we have an exterior product 
\(\ds Y_{1}\mywedge\cdots\mywedge Y_{m}\) of \(\ds Y_{1},\dots,Y_{m}\).  
Omitting $i$-th element, we
have \(\ds Y_{1} \mywedge \cdots \mywedge Y_{i-1} \mywedge Y_{i+1}
\mywedge \cdots \mywedge Y_{m} \), which is often denoted as   
\(\ds Y_{1} \mywedge \cdots\widehat{Y_{i}}\cdots\mywedge Y_{m} \).  

The boundary operator \(\ds \pdel_{} :\myCS{m}\to \myCS{m-1}\)
is defined by 
\begin{align} 
 \pdel ( Y_{1}\mywedge \cdots \mywedge Y_{m}  )
 \kmcomment{
& =  \sum_{i<j} (-1)^{ 
{ \mathop{\sum}_{s<j} (1+ y_{s}y_{j})}  
+ { \mathop{\sum}_{s<i} (1+ y_{s}y_{i})}  } 
 \Sbt{Y_{j}}{Y_{i}} \mywedge \CYm{i,j}{m}   \label{triv:1}
\\
}
&= 
\sum_{i<j} (-1)^{ i-1 + y_{i} \mathop{\sum}_{i< s<j} y_{s}} 
Y_{1} \mywedge \cdots \widehat{ Y_{i} } \cdots \mywedge 
\overbrace{\Sbt{Y_{i}}{Y_{j}}}^{j} \mywedge \cdots  \mywedge Y_{m}
\label{triv:1}
\end{align}
where \(\ds y_{i}\) is the homogeneous degree of 
 \(Y_{i}\), i.e., \(\ds Y_{i} \in \frakg_{y_{i}} \). 
\(\ds
\pdel\circ \pdel = 0\) holds and we have the homology groups 
\[\ds 
\myHom{m}(\frakg, \mR) =  \ker(\pdel : \myCS{m} \rightarrow
\myCS{m-1})/ \pdel ( \myCS{m+1} )\;.  
\] 
\begin{align} 
\noalign{
If all $Y_{i}$ are even (i.e.,  $y_{i}$ are even) in \eqref{triv:1}, then }
 \pdel ( Y_{1}\mywedge \cdots \mywedge Y_{m}  ) 
& = -  \sum_{i<j} (-1)^{ i+j }
 \Sbt{Y_{i}}{Y_{j}} \mywedge 
 Y_{1} \mywedge \cdots\widehat{Y_{i}}\cdots
\widehat{Y_{j}}\cdots \mywedge Y_{m} 
\;. \label{all:even}
\\
\noalign{
If all $Y_{i}$ are odd (i.e., $y_{i}$ are odd) in \eqref{triv:1}, then }
 \pdel ( Y_{1}\mywedge \cdots \mywedge Y_{m}  ) 
& =  \sum_{i<j} 
 \Sbt{Y_{i}}{Y_{j}} \mywedge 
 Y_{1} \mywedge \cdots\widehat{Y_{i}}\cdots
\widehat{Y_{j}}\cdots \mywedge Y_{m}   \;.  \label{all:odd}
\end{align}
\begin{align}
\shortintertext{
Define a binary operation \( \SbtES{A}{B} \) by} 
\SbtES{A}{B} &= 
\pdel ( A\mywedge B ) - (\pdel A)\mywedge B - (-1)^{\bar{a}} A
\mywedge \pdel B 
\label{bdary:bunkai}
\\
\noalign{ If all \(A_{i}\) are even and all \(B_{j}\) are odd, then 
} 
\SbtES{A}{B}
& =   \sum_{i,j} (-1)^{i+1} A_{1}\mywedge\cdots\widehat{ A_{i} }\cdots 
A_{\bar{a}} 
\mywedge \SbtS{A_{i}}{B_{j}} \mywedge
B_{1} \mywedge 
\cdots \widehat{ B_{j} }\cdots \mywedge B_{\bar{b}}\;. 
\end{align}

\section{Superalgebras associated with Lie algebras}

\newcommand{\myXi}[2]{ #1[#2] }
\begin{defn}[Schouten-like bracket]
Let \(\frakg \) be a Lie algebra.  Then we have a \(\mZ\)-graded Lie
superalgebra \(\ovfrakg\) given by \[
\ovfrakg = \sum _{i} \Lambda^{i} \frakg \;,\;\text{ where, the grade of
the element in } \Lambda^{i} \;\text{ is }\; i-1\;,\;
\]
and the bracket
of \(
 a = {a}_{1} \wedge \cdots \wedge {a}_{\bar{a}} \in \Lambda ^{\bar{a}}
\frakg\)  and  \(  
b = {b}_{1} \wedge \cdots \wedge {b}_{\bar{b}} \in
\Lambda ^{\bar{a}} \frakg\) 
is given by 
\begin{align}
\SbtS{a}{b} &= \sum_{i,j} (-1)^{i+j} \Sbt{a_{i}}{b_{j}} \wedge \myXi{a}{i} \wedge
\myXi{b}{j} \\
\noalign{ where }
\myXi{a}{i} &= {a}_{1} \wedge \cdots \wedge {a}_{i-1} \wedge 
{a}_{i+1} \wedge \cdots \wedge  {a}_{\bar{a}} \in \Lambda ^{\bar{a}-1}
\frakg\;\text{ for } i=1,\ldots, \bar{a}\; \\  
\myXi{b}{j} &= {b}_{1} \wedge \cdots \wedge {b}_{j-1} \wedge 
{b}_{j+1} \wedge \cdots \wedge  {b}_{\bar{b}} \in \Lambda ^{\bar{b}-1}
\frakg\; \text{for} j=1,\ldots, \bar{b}\;.   
\end{align}
\end{defn}
\begin{prop}
\begin{equation} \SbtS{b}{a} = - (-1)^{ (\bar{b}-1 )( \bar{a}-1 ) 
} \SbtS{a}{b}
\end{equation}
\begin{align} \SbtS{a}{b\wedge c}  =&  \SbtS{a}{b} \wedge c  + 
 (-1)^{ (\bar{a}-1 )\bar{b}  } b \wedge \SbtS{a}{c} \label{eqn:deri:1}  \\
\SbtS{a\wedge b}{c}  =&  a \wedge \SbtS{b}{c}   + 
 (-1)^{ \bar{b} (\bar{c} -1)  } \SbtS{a}{c}   \wedge b \label{eqn:deri:2}
\end{align}
\begin{equation}
\mathop{\mathfrak S}_{a,b,c} (-1)^{(\bar{a}-1)(\bar{c}-1)} 
\SbtS{\SbtS{a}{b}}{c}= 0\end{equation}
\end{prop}

\section{Homology groups of superalgebras associated with Lie algebras}
Let \(\frakg\) be a finite dimensional Lie algebra with the Lie bracket
\([\cdot,\cdot]\).  Then we have chain complex consisting of \(\CSp{i}=
\Lambda^{i}\frakg\) and homology groups. 

\(\ds\ovfrakg = \mathfrak{\oplus}_{i=0}^{-1 +\dim \frakg} \frakg_{i} \) becomes a \(\mZ\)-graded
Lie superalgebra (where \(\frakg_{i} = \CSp{i+1}\)) 
with the Schouten bracket \(\SbtS{\cdot}{\cdot}\).
From the definition, $w$-weighted  $m$-th chain space  
\(\wtedC{m}{w}\)
is defined as follows:
\begin{align}
 \wtedC{m}{w} &= \sum \mywedge^{p_{0}} \frakg_{0}
 \mywedge^{p_{1}} \frakg_{1}
\mywedge^{p_{2}} \frakg_{2}\mywedge \cdots 
\mywedge^{p_{s}} \frakg_{s} \text{ where } s = -1+ \dim \frakg 
\notag
\\
m &= p_{0} + p_{1} + \cdots + p_{s}\;, 
\quad
w = 0 p_{0} + 1 p_{1} + \cdots + s p_{s}\;, 
\notag
\\ p_{2i} &\leqq \binom{\dim\frakg}{2i+1}\text{ , which is the dimensional
 restriction for } \frakg_{2 i}\;.  
\label{dim:cond}\\\noalign{In particular}
\wtedC{m}{0} &= \mywedge^{m} \frakg_{0} \;\text{ for }\; m=0,\ldots,
\dim \frakg\; . 
\end{align}
By the same discussion we already developed, we see 
\[ w + m = 1 p_{0} + 2 p_{1} + \cdots + ({\dim\frakg}) p_{s} \]
and we deal with the Young diagrams of the area \(w + m\) with the
length \(m\) and pick up those satisfy the dimensional condition
\eqref{dim:cond}.  
Since 
\(p_{0}\) does not contribute for the weight \(w\), we assume
\(p_{0}=0\), and for given \(w\) we define the subspace  
\(\ds \wtedCR{m}{w}\) of  \(\ds \wtedC{m}{w}\) by 
\[ 
 \wtedCR{m}{w} = \sum
 \mywedge^{p_{1}} \frakg_{1}
\mywedge^{p_{2}} \frakg_{2}\mywedge \cdots 
\mywedge^{p_{s}} \frakg_{s} \text{ where }
 \sum_{i}i p_{i} = w\text{ , }  \sum_{i} p_{i} = m\;,
\text{  and satisfy \eqref{dim:cond}.}  \] 

In order to control 
\(\ds \wtedCR{m}{w}\), we define the subspace  
\(\ds \wtedCRx{m}{w}\) of \(\ds \wtedC{m}{w}\) by 
\(p_{0} = p_{1}=0\). Thus,  
\[ 
 \wtedCRx{m}{w} = \sum
\mywedge^{p_{2}} \frakg_{2}\mywedge \cdots 
\mywedge^{p_{s}} \frakg_{s} \text{ where }
 \sum_{i>1}i p_{i} = w\text{ , }  \sum_{i>1} p_{i} = m\;,
\text{  and satisfy \eqref{dim:cond}.}  \] 
We define subspaces \(\ds \wtedCR{\bullet}{w} = \sum_{m} \wtedCR{m}{w} \)  
and \(\ds \wtedCRx{\bullet}{w} = \sum_{m} \wtedCRx{m}{w} \).   
Then we have
\begin{align}
\wtedC{\bullet}{w} & = 
\wtedC{\bullet}{0} \mywedge  
\wtedCR{\bullet}{w}\;,\quad\text{ precisely }
 \wtedC{m}{w} = 
 \sum_{a+b=m} 
 \mywedge^{a} \frakg_{0} \mywedge 
 \wtedCR{b}{w}\; , \label{eqn:g0:others}\\
\wtedCR{\bullet}{w} & = 
\wtedCRx{\bullet}{w} +
\frakg_{1} \mywedge \wtedCR{\bullet}{w-1}   
\;, \quad   
\wtedCRx{\bullet}{1}  = 0\;,\;
\wtedCR{\bullet}{0}  = 1\;,\; 
\wtedCR{\bullet}{1} =  \frakg_{1}\; . \label{eqn:rist:recursive} 
\end{align}
\eqref{eqn:g0:others} implies the next theorem:  
\begin{thm}
The Euler number of \(w\)-weighted homology groups of superalgebra
induced from finite dimensional Lie algebra is 0.
\end{thm}
\textbf{Proof:} 
 \begin{align*}
 \text{The Euler num} &= \sum_{m} (-1)^{m} \dim 
 \wtedC{m}{w} = 
 \sum_{a+b=m} (-1)^{a+b} \dim(\mywedge^{a} \frakg_{0}) \dim \wtedCR{b}{w}  
\\&= 
 \sum_{a} (-1)^{a} \dim(\mywedge^{a} \frakg_{0}) \sum_{b} (-1)^{b} \dim \wtedCR{b}{w}  
\\&= 
 \sum_{a} (-1)^{a} \binom{\dim \frakg_{0}}{a} \sum_{b} (-1)^{b} \dim
 \wtedCR{b}{w} = 0\;.  
 \end{align*} \kmqed

The equation 
\eqref{eqn:rist:recursive} shows 
the space \( \wtedCR{\bullet}{w} \) is determined by 
\( \wtedCRx{\bullet}{w} \) and 
\( \wtedCR{\bullet}{w-1} \), and suggests discussions by induction on the
weight \(w\).

\section{Super homologies of 2-dimensional Lie algebras}
Since our discussion is trivial for abelian Lie algebras, we consider a Lie
algebra \(\frakg\)  
generated by \( \zb{1}, \zb{2} \) with \( \Sbt{ \zb{1} }{\zb{2}} =
\zb{1}\). We have only 1-dimensional 2-vectors, so we take \(\ub{}= \zb{1}
\wedge \zb{2}\).  
For a given weight \(w\), the chain complex becomes as follows. 
Using 
\( \pdel ( \zb{1} \mywedge \zb{2}  ) = \zb{1}\), 
\( \SbtES{ \zb{1} }{ \ub{} } = 0\) and 
\( \SbtES{ \zb{2} }{ \ub{} } = - \ub{}\),  we see 
\( \SbtES{ \zb{1} }{ \ub{}^{w} } = 0\) and 
\( \SbtES{ \zb{2} }{ \ub{}^{w} } = - w \ub{}^{w}\).  We have \( \pdel \ub{}
^{w} =0\),  
\( \pdel ( \zb{1} \mywedge \ub{}^{w} ) = \SbtES{\zb{1}}{ \ub{} ^{w} } = 0\), 
\( \pdel ( \zb{2} \mywedge \ub{}^{w} ) = \SbtES{\zb{2}}{ \ub{} ^{w} }  
= - w \ub{} ^{w} \), and 
 \( \pdel ( \zb{1} \mywedge \zb{2} \mywedge  \ub{}^{w})  = (1+w) \zb{1}
 \mywedge \ub{}^{w}\). Using those, we complete the table below.  

\begin{center}
\(
\begin{array}{c | *{5}{c}}
 &  \wtedC{0}{0} & \wtedC{1}{0} & \wtedC{2}{0} \\\hline\hline
 \text{generators}  & 1 & \zb{i} 
& \zb{1}\mywedge \zb{2}  \\ \hline
\text{SpaceDim} & 1 & 2 & 1 \\\hline
\text{KerDim} & 1 & 2 & 0 \\\hline
\text{Betti} & 1 & 1 & 0 \\\hline
 \end{array}
\)
\hfil
\(
\begin{array}{c | *{5}{c}}
w>0  &  \wtedC{w}{w} & \wtedC{w+1}{w} & \wtedC{w+2}{w} \\\hline\hline
 \text{generators}  & \ub{}^{w} & \zb{i}\mywedge \ub{}^{w} 
& \zb{1}\mywedge \zb{2} \mywedge \ub{}^{w} \\ \hline
\text{SpaceDim} & 1 & 2 & 1 \\\hline
\text{KerDim} & 1 & 
1 
 & 0 \\\hline
\text{Betti} & 0 & 0 & 0 \\\hline
 \end{array}
\)

\end{center}

\section{Super homologies of 3-dimensional Lie algebras}
We know  the all Lie algebras of lower dimensional cases (cf.\ Lie algebras
by Jacobson), in this section we try to get super homology tables for
3-dimensional Lie algebras. For that purpose, we prepare notations here. 
Let  \(\zb{1}, \zb{2}, \zb{3}\) be a basis of \(\frakg\),   
\(\ub{1}, \; \ub{2},\; \ub{3}\) be 
a basis of \(\Lambda^{2}\frakg\), and \(\zb{4} = V = \zb{1} \wedge 
 \zb{2} \wedge  \zb{3} \) be a basis of \(\Lambda^{3}\frakg\). 
Since \(\{\ub{i}\}\) have symmetric property, 
we put \(U^{a,b,c} = \ub{1}^{a}\mywedge
\ub{2}^{b}\mywedge  \ub{3}^{c}\). 
\kmcomment{
\(\mywedge\), we put \(U^{a,b,c} = \mywedge^{a}\ub{1}^{a}\mywedge
\mywedge^{b}\ub{2}^{b}\mywedge  \mywedge^{c}\ub{3}^{c}\). 
}

Denote the 3-vector \((a,b,c)\) by \(A\) and also \(\myUnital{j}\) means
3-vector \( (\delta^{k}_{j})_{k=1..3}\), and \( A - \myUnital{j}\) means
the difference of 3-vectors. 
To handle even basis at once, we introduce
an odd notation:
\begin{equation}
 W^{\eps{1},\eps{2},\eps{3},\eps{4}}
= \zb{1}^{\eps{1}}\mywedge \zb{2}^{\eps{2}}\mywedge
\zb{3}^{\eps{3}}\mywedge \zb{4}^{\eps{4}} 
\label{odd:express:new}
\end{equation}
where \(\eps{1},\eps{2},\eps{3},\eps{4}\) are 0 or 1, and if
\(\eps{1}=1\) then \(\zb{1}^{\eps{1}}= \zb{1}\)
as expected otherwise we require \( \zb{1}^{\eps{1}}\) disappear or becomes
nothing/``empty''.  For instance, \( W^{1,0,1,0} = \zb{1}
\mywedge \zb{3}\). 
A generic chain of the chain complex of superalgebra is 
\begin{equation} 
W^{\eps{1},\eps{2},\eps{3},\eps{4}} \mywedge \bU{a,b,c}{\ell}
\label{gen:form}
\end{equation}
The degree(length) is 
\(\sum_{i=1}^{4}\eps{i} + a+b+c \) and  the weight is 
 \(2\eps{4} + a+b+c \) for the chain \eqref{gen:form}.   

\renewcommand{\arraystretch}{0.8}
For a given weight \(w\), the chain complex becomes as follows: 
\begin{center}
\(
\begin{array}{c | *{5}{c}}
 & \wtedC{w-1}{w} & \wtedC{w}{w} & \wtedC{w+1}{w} & \wtedC{w+2}{w} &
 \wtedC{w+3}{w} 
\\\hline\hline
 W^{\varE} \mywedge \bU{A}{w} & \text{none} & 
 \begin{tabular}{c}
 \( \epsilon\epsilon=0\)\\ \( \eps{4}=0\) \end{tabular} &
 \begin{tabular}{c}
 \( \epsilon\epsilon=1\)\\ \( \eps{4}=0\) \end{tabular} &
 \begin{tabular}{c}
 \( \epsilon\epsilon=2\)\\ \( \eps{4}=0\) \end{tabular} &
 \begin{tabular}{c}
 \( \epsilon\epsilon=3\)\\ \( \eps{4}=0\) \end{tabular} 
 \\
 \hline
 W^{\varE} \mywedge \bU{P}{w-2} & 
 \begin{tabular}{c}
 \( \epsilon\epsilon=0\)\\ \( \eps{4}=1\) \end{tabular} &
 \begin{tabular}{c} \( \epsilon\epsilon=1\)\\ \( \eps{4}=1\) \end{tabular} 
 &
 \begin{tabular}{c}
 \( \epsilon\epsilon=2\)\\ \( \eps{4}=1\) \end{tabular} &
 \begin{tabular}{c}
 \( \epsilon\epsilon=3\)\\ \( \eps{4}=1\) \end{tabular} &
 \text{none}  
\\\hline\hline
 W^{\varE} \mywedge \bU{A}{w} & \text{none} & 
 \bU{A}{w}&
 \zb{i}\mywedge \bU{A}{w}&
 \zb{i}\mywedge \zb{j}\mywedge \bU{A}{w}&
W^{1110}
 \mywedge  \bU{A}{w} 
 \\ \hline
 W^{\varE} \mywedge \bU{P}{w-2} & 
 \zb{4} \mywedge \bU{P}{w-2}
 &
 \zb{i} \mywedge \zb{4} \mywedge \bU{P}{w-2}
 &
 W^{\eps{i}=0}\mywedge \bU{P}{w-2}
 &
 W^{1111} \mywedge  \bU{P}{w-2}
 & \text{none}  \\ \hline\hline
 \text{SpaceDim} & \binom{w}{2} & 3 \binom{w}{2} +  \binom{w+2}{2}  
 & 3 \binom{w}{2} + 3 \binom{w+2}{2}  
 &  \binom{w}{2} + 3 \binom{w+2}{2}  
 &   \binom{w+2}{2}  \\\hline
 \end{array}
\)
\end{center}
where 
\(\varE = (\eps{1},\eps{2},\eps{3},\eps{4})\) and  
\(\epsilon\epsilon=\sum_{i=1}^{3} \eps{i}\).
\renewcommand{\arraystretch}{1.0}
In general, the boundary operator works as below:  
\begin{align}
\pdel ( W^{\varE} \wedge U^{A} ) &= 
\pdel ( W^{\varE}) \wedge U^{A}  
+ (-1)^{|\varE|}  W^{\varE} \wedge \pdel  U^{A}  
+ \SbtES{ W^{\varE}}{  U^{A}  } 
\\
\pdel  W^{\varE} &= - \sum _{i<j} \myeps{i} \myeps{j} \SbtS{\zb{i}}{\zb{j}}
\mywedge W ^{ \varE - \myUnital{i} - \myUnital{j} }\;, \quad 
\text{
here } \myeps{i} = \eps{i} (-1)^{\sum_{k\leqq i} \eps{k}}  
\\
 \SbtES{ W^{\varE}}{U^{A}} &= 
 - \sum_{i\leqq 3} \myeps{i} W^{\eps{i}=0} \mywedge \SbtES{\zb{i}}{ U^{A} }
\\
\pdel U^{A} &= 
\sum_{i} \binom{a_{i}}{2} \SbtS{\ub{i}}{\ub{i}} \mywedge U^{A
- 2 \myUnital{i}}
+
\sum_{i<j} {a_{i}}{a_{j}} \SbtS{\ub{i}}{\ub{j}} \mywedge U^{A
-  \myUnital{i} -  \myUnital{j} }
\end{align}

\paragraph{ The boundary operator on \(\wtedC{w-1}{w}\) is trivial:} 
\null

A general chain is given by 
\( M = \sum_{P}  
\MC{}{P} \zb{4} \mywedge \bU{P}{w-2} \) and  
\( \pdel M = - \sum_{P} \MC{}{P} \zb{4} \mywedge \pdel \bU{P}{w-2} 
+ \sum_{P} \MC{}{P}  \SbtES{\zb{4}} {\bU{P}{w-2}} = 0\) 
because \( \pdel \bU{P}{w-2} = \zb{4}\mywedge \text{"some"}\) and \(
\SbtS{\zb{4}}{ \ub{i} } = 0\). Thus, we have
\begin{prop}
\label{prop:triv}
Even though the boundary operator is really depending on the Schouten
bracket,  the boundary operator on 
the lowest chain space 
is trivial. 
\end{prop}
 
\paragraph{ The boundary operator on \(\wtedC{w}{w}\):} \null  
\begin{align*}
&\quad \pdel( \sum \LC{} {A} \bU{A}{w}+ \sum \MC{j}{P} \zb{i} \mywedge \zb{4}
\mywedge \bU{P}{w-2}) 
\\&= 
\sum \LC{} {A} \pdel \bU{A}{w} 
+ \sum \MC{i}{P} ( 
\SbtS{\zb{i}}{\zb{4}} \mywedge \bU{P}{w-2}
+ \zb{i} \mywedge \zb{4} \mywedge \pdel \bU{P}{w-2} 
+ \SbtES{\zb{i} \mywedge \zb{4}} {\bU{P}{w-2}}
) 
\\&= 
\sum \LC{} {A} \pdel \bU{A}{w} 
+ \sum \MC{i}{P}  
\SbtS{\zb{i}}{\zb{4}} \mywedge \bU{P}{w-2}
 + \zb{4} \mywedge \sum \MC{i}{P}  \SbtES{\zb{i}} {\bU{P}{w-2}}
\end{align*}
\paragraph{ The boundary operator on \(\wtedC{w+1}{w}\):} 
\begin{align*}
&\quad \pdel( \sum \LC{i} {A}\zb{i} \mywedge \bU{A}{w}
+ \sum \MC{j}{P} W^{\eps{j}=0} \mywedge \bU{P}{w-2}) 
\\&= 
\sum \LC{i} {A} ( -\zb{i} \mywedge  \pdel \bU{A}{w} + \SbtES{\zb{i}}{\bU{A}
{w}}) 
+ \sum \MC{j}{P} ( \pdel W^{\eps{j}=0}) \mywedge \bU{P}{w-2} + \sum \MC{j}{P} 
 \SbtES{
 W^{\eps{j}=0}}{\bU{P}{w-2}}
\end{align*}

\paragraph{ The boundary operator on \(\wtedC{w+2}{w}\):} 
\begin{align*}
&\quad \pdel( \sum \LC{i} {A} W^{\eps{i}=\eps{4}=0} \mywedge \bU{A}{w}
+ \sum \MC{ }{P} W^{1111} \mywedge \bU{P}{w-2}) 
\\&= 
\sum \LC{i}{A}  
(\pdel W^{\eps{i}=\eps{4}=0})\mywedge \bU{A}{w}
+ \sum \LC{i}{A}  
  W^{\eps{i}=\eps{4}=0} \mywedge \pdel \bU{A}{w}
+
\sum \LC{i}{A}  
  \SbtES{ W^{\eps{i}=\eps{4}=0}}{\bU{A}{w}} 
  \\& \quad 
+ \sum \MC{ }{P} \pdel (W^{1111}) \mywedge \bU{P}{w-2} 
+ \sum \MC{ }{P} \SbtES{ W^{1111} }{\bU{P}{w-2}}  
\end{align*}

\paragraph{ The boundary operator on \(\wtedC{w+3}{w}\):} 
\begin{align*}
&\quad \pdel  \sum \LC{ } {A} W^{\eps{4}=0} \mywedge \bU{A}{w}
\\&= 
\sum \LC{i}{A}  
(\pdel W^{\eps{4}=0})\mywedge \bU{A}{w}
- \sum \LC{i}{A}  W^{\eps{4}=0} \mywedge \pdel \bU{A}{w}
+ \sum \LC{i}{A}  \SbtES{ W^{\eps{4}=0}}{\bU{A}{w}} 
\\&=  
( 
\SbtS{ \zb{1} }{\zb{2}} \mywedge \zb{3} 
-\SbtS{ \zb{1} }{\zb{3}} \mywedge \zb{2} 
+\SbtS{ \zb{2} }{\zb{3}} \mywedge \zb{1} 
)\mywedge \sum \LC{i}{A}  \bU{A}{w}
- \sum  W^{\eps{4}=0} \mywedge \LC{i}{A} \pdel \bU{A}{w}
\\&\quad 
+ W^{0110} \mywedge \LC{}{A} \SbtES{\zb{1}}{\bU{A}{w}}
- W^{1010} \mywedge \LC{}{A} \SbtES{\zb{2}}{\bU{A}{w}}
+ W^{1100} \mywedge \LC{}{A} \SbtES{\zb{3}}{\bU{A}{w}}
\end{align*}

\subsection{\(\dim\frakg=3\), \(\dim [\frakg,\frakg]=1\) and \( [\frakg,\frakg] \subset Z(\frakg)\) }

Consider a Lie algebra \(\frakg\) where \([\frakg,\frakg]\) is 1-dimensional
and is in of the center of \(\frakg\). Then we find a basis \(\zb{1}, \zb{2},
\zb{3}\) of \(\frakg\) so that \( \Sbt{\zb{1}}{\zb{2}} = -\Sbt{\zb{2}}
{\zb{1}} = \zb{3}\) and the other brackets are zero.  We take \(\ub{1}=
\zb{2}\wedge\zb{3},\; \ub{2}= \zb{3}\wedge\zb{1},\; \ub{3}= \zb{1}
\wedge\zb{2}\) as a basis of \(\Lambda^{2}\frakg\).  It should be noted that
this definition does not necessarily follow a natural order. 
Take \(V = \zb{4} = \zb{1} \wedge 
 \zb{2} \wedge  \zb{3} \) as a basis of \(\Lambda^{3}\frakg\). 
Now we have the multiplication (by the Schouten bracket) tables:
\begin{center}
\( \begin{array}[t]{c|*{3}{c}|c|*{3}{c}|}
       & \zb{1} & \zb{2} & \zb{3}&\zb{4} & \ub{1}&\ub{2}& \ub{3}  \\\hline 
\zb{1} &  0     & \zb{3} & 0 & 0     & 0     & 0    & -\ub{2} \\
\zb{2} & -\zb{3}& 0      & 0 & 0     & 0     & 0    & \ub{1}  \\
\zb{3} & 0 & 0 & 0  & 0 & 0 & 0 & 0  
\end{array}
 \)
\hfil
\( \begin{array}[t]{c|*{3}{c}|c}
 & \ub{1} & \ub{2} & \ub{3} & \zb{4} \\\hline 
\ub{1} & 0 & 0 & 0 & 0 \\
\ub{2} & 0 & 0 & 0 & 0 \\
\ub{3} & 0 & 0 & 2\zb{4}  & 0\\
\end{array}
 \)
\end{center} 
Then we have 
\begin{align*}
\pdel W^{1111} &= 0\;,\quad 
\pdel W^{1110} = 0 \;,\quad 
\pdel W^{1101} = \zb{3} \mywedge \zb{4} \;,\quad 
\pdel W^{1011} =  0 \;,\quad 
\pdel W^{0111} =  0 \;,\quad 
\\
\pdel U^{a,b,c} &= \zb{4} \mywedge 
  2 \tbinom{c}{2}  U^{a, b, c-2} 
\;. 
\quad 
\sum_{A} \LC{}{A}\pdel U^{a,b,c} = \zb{4} \mywedge \sum_{} U^{p,q,r} 
  (r+2) (r+1) \LC{}{p, q, r+2} 
\;. 
\\
\SbtES{\zb{1}}{U^{A}} &=  - c U^{a,b+1,c-1}\;,  \quad  
\SbtES{\zb{2}}{U^{A}}  = c U^{a+1,b,c-1}\;,  \quad 
\SbtES{\zb{3}}{U^{A}} = 0 
\\ \sum_{A} \LC{i}{A} \SbtES{\zb{1}}{U^{A}} 
&= \sum_{A}(-(c+1)\LC{i}{a,b-1,c+1}) U^{A}\;,\quad 
 \sum_{A} \LC{i}{A} \SbtES{\zb{2}}{U^{A}} 
= \sum_{A}( (c+1)\LC{i}{a-1,b,c+1}) U^{A}\;.
\end{align*}

\paragraph{Kernel of 
\(\ds\wtedC{w+3}{w}\mathop{\to}^{\pdel} \wtedC{w+2}{w} \)}
A basis of \( \wtedC{w+3}{w} \), we have 
\(W^{1,1,1,0 } \mywedge U^{a,b,c} \) with \( a+b+c= w\) and  
take a linear combination
\( X =  \sum_{a+b+c=w} \LC{}{a,b,c} 
W^{1,1,1,0 } \mywedge U^{a,b,c} \) by unknown scalars  \(\LC{}{a,b,c}
\). 
Then  
\begin{align*} \pdel X =&  \sum_{a+b+c= w} 
\LC{}{a,b,c} \left(
c W^{0,1,1,0 } \mywedge U^{a,b+1,c-1} 
+ c  W^{1,0,1,0 } \mywedge U^{a+1,b,c-1} \right.
\\&\qquad \qquad \qquad 
\left. 
- c (c-1)
W^{1,1,1,1 } \mywedge U^{a,b,c-2} 
\right)
\\
= &  
 W^{0,1,1,0 }  \sum_{a+b+c= w} \LC{}{a,b,c}
c \mywedge U^{a,b+1,c-1} 
+  W^{1,0,1,0 }  \sum_{a+b+c= w} \LC{}{a,b,c} c \mywedge U^{a+1,b,c-1}
\\&\quad 
- W^{1,1,1,1 }   \sum_{a+b+c= w} \LC{}{a,b,c}
 c (c-1)  \mywedge U^{a,b,c-2} 
\end{align*}
\(\pdel X = 0\) is generated by 
\begin{subequations}
\begin{align} \label{G3D1Y:A:1}
& \LC{}{a,b-1,c+1} (c+1)   \\
\label{G3D1Y:A:2}
& \LC{}{a-1,b,c+1} (c+1)  \\
\label{G3D1Y:A:3}
& \LC{}{p,q,r+2} (r+2) (r+1)  
\end{align}
\end{subequations}
\kmcomment{
\eqref{G3D1Y:A:1} or 
\eqref{G3D1Y:A:2} show that \(\LC{}{a,b,c} =0 \) if \( a+b \leq w-1\), and 
\eqref{G3D1Y:A:3} shows   \(\LC{}{a,b,c} =0 \) if \( a+b \leq w-2\).  Thus, 
 \(\LC{}{a,b,c} =0 \) if \( a+b \leq w-1\), and 
 \(\LC{}{a,b,c}\) are free  if \( a+b =w \), and  the kernel dimension
 is \(w+1\). 

Revise on July 20, 2020 by focusing the rank, the number of generators of
linear equations: 

\eqref{G3D1Y:A:1} $\sim$ 
\eqref{G3D1Y:A:3}  yields the generators are
\[
 (c+1)\lambda_{a,b-1,c+1}  \;,\; 
 (c+1)\lambda_{a-1,b,c+1}  \;,\;
  (r+2)(r+1)  \lambda_{p,q,r+2}
\]
}
Thus, the rank is \(\binom{w+2}{2} - \# \{ \LC{}{a,b,0} \mid a+b=w \} =
\binom{w+2}{2} - (w+1)\) 
from \eqref{G3D1Y:A:1} $\sim$ \eqref{G3D1Y:A:3},   
and the kernel dimension is \(w+1\).

\paragraph{Kernel of \(\ds\wtedC{w+2}{w}\mathop{\to}^{\pdel} \wtedC{w+1}{w} \)}
As a basis of \( \wtedC{w+2}{w} \), we have 
\(
W^{ \eps{1}, \eps{2} ,\eps{3},0} \mywedge \bU{a,b,c}{w}\) with 
\( \eps{1}+ \eps{2} +\eps{3}= 2\) 
and 
\(W^{ 1,1,1 ,1} \mywedge \bU{p,q,r}{w-2}\).   
\kmcomment{
手原稿の \labda^{1} -- [1,1,0] --> \newlabmda^{3}, 
手原稿の \labda^{2} -- [1,0,1] --> \newlabmda^{2}, 
手原稿の \labda^{3} -- [0,1,1] --> \newlabmda^{1}, 
と変更する。1 と 3 の交換のみ。
}

Take a linear combination
\[ X =  \sum_{a+b+c=w} \LC{k}{a,b,c} 
W^{ \eps{1}, \eps{2} ,\eps{3},0} \mywedge \bU{a,b,c }{w}
+  \sum_{p+q+r=w-2} \MC{}{p,q,r} W^{ 1,1,1,1} \mywedge \bU{p,q,r}{w-2}
\] 
by unknown scalars  \(\lambda^{k}_{a,b,c},  \mu_{p,q,r}\) 
with 
\( \eps{1}+ \eps{2} +\eps{3}= 2\) and \(\eps{k} = 0\) .
\begin{align*}
\pdel X =&  
\sum \LC{1}{a,b,c} ( 
c W^{0,0,1,0} \bU{a+1,b,c-1}{w}
+ 2 \binom{c}{2}  W^{0,1,1,1} \bU{a,b,c-2}{w-2})
\\&
+ \sum \LC{2}{a,b,c} ( -
c W^{0,0,1,0} \bU{a,b+1,c-1}{w}
+2\binom{c}{2}  W^{1,0,1,1} \bU{a,b,c-2}{w-2})
\\&
+ \sum \LC{3}{a,b,c} ( W^{0,0,1,0} \bU{a,b,c}{w}
- c W^{0,1,0,0} \bU{a,b+1,c-1}{w}
- c W^{1,0,0,0} \bU{a+1,b,c-1}{w}
+ 2 \binom{c}{2} W^{1,1,0,1} \bU{a,b,c-2}{w-2})
\\& 
-
\sum \MC{}{p,q,r} ( 
r W^{0,1,1,1} \bU{p,q+1,r-1}{w-2}
+ r W^{1,0,1,1} \bU{p+1,q,r-1}{w-2}
)
\\ =&  
   W^{0,0,1,0}(\sum \LC{1}{a,b,c}  c \bU{a+1,b,c-1}{w}
-   \sum \LC{2}{a,b,c}  c \bU{a,b+1,c-1}{w}
+   \sum \LC{3}{a,b,c} \bU{a,b,c}{w})
\\& 
-  W^{0,1,0,0} \sum \LC{3}{a,b,c} c \bU{a,b+1,c-1}{w}
-  W^{1,0,0,0} \sum \LC{3}{a,b,c} c \bU{a+1,b,c-1}{w}
\\& 
+  W^{0,1,1,1} (\sum \LC{1}{a,b,c}  2 \binom{c}{2} \bU{a,b,c-2}{w-2}
-   \sum \MC{}{p,q,r}  r \bU{p,q+1,r-1}{w-2})
\\& 
+  W^{1,0,1,1} (\sum \LC{2}{a,b,c}  2\binom{c}{2} \bU{a,b,c-2}{w-2}
-   \sum \MC{}{p,q,r}  r \bU{p+1,q,r-1}{w-2})
+  W^{1,1,0,1} \sum \LC{3}{a,b,c} 2 \binom{c}{2} \bU{a,b,c-2}{w-2}
\end{align*}
Thus, 
 \( \pdel X=0\) is determined by 

\begin{subequations}
\begin{align} \label{G3D1Y:Rev:1}
& \LC{1}{a-1,b,c+1}  (c+1) -   \LC{2}{a,b-1,c+1}  (c+1) + \LC{3}{a,b,c}
\\& \LC{3}{a,b-1,c+1} (c+1)  \label{G3D1Y:Rev:2}
\\& \LC{3}{a-1,b,c+1} (c+1)  \label{G3D1Y:Rev:3}
\\& \LC{1}{p,q,r+2} 2 \binom{r+2}{2} -  \MC{}{p,q-1,r+1}  (r+1)  \label{G3D1Y:Rev:4}
\\& \LC{2}{p,q,r+2} 2 \binom{r+2}{2} -  \MC{}{p-1,q,r+1}  (r+1 )  \label{G3D1Y:Rev:5}
\\& \LC{3}{p,q,r+2} 2 \binom{r+2}{2} \;.  \label{G3D1Y:Rev:6}
\end{align}
\end{subequations}
\eqref{G3D1Y:Rev:2}, 
\eqref{G3D1Y:Rev:3}, and  
\eqref{G3D1Y:Rev:6} tell that  \(\LC{3}{a,b,c} \) 
are single generators for \( c>0\). With \( c=0\) for \eqref{G3D1Y:Rev:1},
\(\LC{3}{A} \) are linear independent generators.  
When \(c > 0\) for \eqref{G3D1Y:Rev:1}, then we have 
\[
 \LC{1}{a-1,b,c+1}   -   \LC{2}{a,b-1,c+1}  
\]
Those are obtained from \eqref{G3D1Y:Rev:4} and \eqref{G3D1Y:Rev:5}. Thus, the rank is 
\( \binom{w+2}{2} + 2 \binom{w}{2} \), the kernel dimension is  
\( ( \binom{w}{2}+ 3 \binom{w+2}{2}) - (\binom{w+2}{2} + 2 \binom{w}{2}) = 
-  \binom{w}{2} + 2 \binom{w+2}{2}  \).

\paragraph{Kernel of 
\(\ds\wtedC{w+1}{w}\mathop{\to}^{\pdel} \wtedC{w}{w} \)}
As a basis of \( \wtedC{w+1}{w} \), we have 
\( \zb{i} \mywedge \bU{a,b,c}{w} \) with \(i=1,2,3\),  and 
\(W^{ \eps{1},\eps{2},\eps{3} ,1} \mywedge \bU{p,q,r}{w-2}\) with   
\( \eps{1}+ \eps{2} +\eps{3}= 2\).  
Take a linear combination 
\[ X = \sum \LC{i}{a,b,c} \zb{i} \mywedge \bU{a,b,c}{w} + 
 \sum_{ \eps{j}=0, \eps{1}+ \eps{2} +\eps{3}= 2 }
 \MC{j}{p,q,r} 
W^{ \eps{1},\eps{2},\eps{3} ,1} \mywedge \bU{p,q,r}{w-2}\]
for unknown scalars 
\(\LC{i}{a,b,c}\) and \(\MC{j}{p,q,r}\).  
\begin{align*}
\pdel X =& 
\sum \LC{1}{a,b,c} \left(
c\bU{a,b+1,c-1}{w} - 2\binom{c}{2} \zb{1}\mywedge V \mywedge \bU{a,b,c-2}{w-2} 
\right)
- \sum \LC{2}{a,b,c} \left(
c\bU{a+1,b,c-1}{w} + 2\binom{c}{2} \zb{2}\mywedge V \mywedge \bU{a,b,c-2}{w-2} 
\right)
\\ & 
- \sum \LC{3}{a,b,c}  \left( 
2 \binom{c}{2}   \zb{3}\mywedge V \mywedge \bU{a,b,c-2}{w-2} 
\right)
- \sum \MC{1}{p,q,r} \left( 
  r \zb{3}\mywedge V \mywedge \bU{p+1,q,r-1}{w-2} 
\right)
\\ & 
+ \sum \MC{2}{p,q,r} \left( 
  r  \zb{3}\mywedge V \mywedge \bU{p,q+1,r-1}{w-2} 
\right)
 + \sum \MC{3}{p,q,r} 0 \quad\text{(remember that 
\( \MC{3}{p,q,r} \) disappeared.) }
\\ =& 
\sum \LC{1}{a,b,c} c\bU{a,b+1,c-1}{w} - \sum \LC{2}{a,b,c} c\bU{a+1,b,c-1}{w} 
- \zb{1}\mywedge V \mywedge \sum \LC{1}{a,b,c} 2\binom{c}{2} \bU{a,b,c-2}{w-2} 
\\ & 
- \zb{2}\mywedge V \mywedge \sum \LC{2}{a,b,c} 2\binom{c}{2} \bU{a,b,c-2}{w-2} 
\\ & 
- \zb{3}\mywedge V \mywedge 
\left(
\sum \MC{1}{p,q,r}  r \bU{p+1,q,r-1}{w-2} 
-  \sum \MC{2}{p,q}  r  \bU{p,q+1,r-1}{w-2} 
+  \sum \LC{3}{a,b,c} 2 \binom{c}{2}\bU{a,b,c-2}{w-2} 
\right) 
\end{align*}
Thus, we have linear equations 
\begin{align}
& 
\sum \LC{1}{a,b,c} c\bU{a,b+1,c-1}{w} - \sum \LC{2}{a,b,c} c\bU{a+1,b,c-1}{w} = 0 \tag{F1} \label{G3D1Y:NL:E1}
\\& \sum \LC{1}{a,b,c} 2\binom{c}{2} \bU{a,b,c-2}{w-2} = 0 \tag{F2} \label{G3D1Y:NL:E2}
\\& \sum \LC{2}{a,b,c} 2\binom{c}{2} \bU{a,b,c-2}{w-2} =0 \tag{F3} \label{G3D1Y:NL:E3}
\\& 
\sum \MC{1}{p,q,r}  r \bU{p+1,q,r-1}{w-2} 
-  \sum \MC{2}{p,q,r}  r  \bU{p,q+1,r-1}{w-2} 
+  \sum \LC{3}{a,b,c} 2 \binom{c}{2}\bU{a,b,c-2}{w-2} 
= 0 \tag{FF} \label{G3D1Y:NL:EF}
\end{align}
\kmcomment{
\eqref{G3D1Y:NL:E2} and 
\eqref{G3D1Y:NL:E3} imply \begin{equation}
\LC{1}{a,b} = 0\;,\;  
\LC{2}{a,b} = 0\quad\text{ for }\   a+b \leqq w-2 \;. \label{G3D1Y:join:F2:F3}
\end{equation}
Now, 
\eqref{G3D1Y:NL:E1} implies 
\( \LC{1}{a,b} \) and  
\( \LC{2}{a,b} \) with \(a+b=w\) have no restriction, namely they are
free. And also  \eqref{G3D1Y:join:F2:F3} and 
\eqref{G3D1Y:NL:E1} implies 
\begin{equation}
\sum_{a+b=w-1} \LC{1}{a,b} \bU{a,b+1}{w} - \sum_{a+b=w-1} \LC{2}{a,b} \bU{a+1,b}{w} = 0 
\end{equation} and we see that
\[
\LC{1}{0,w-1} = 0\;,\; 
\LC{1}{a,w-1-a} = \LC{2}{a-1,w-a}\; (0<a<w-1)\;, \; \LC{2}{w-1,0}= 0 \; 
\]
Thus, the number of free parameters among 
\( \LC{1}{a,b} \) and  \( \LC{2}{a,b} \) is  \( 2 (w+1) + (w-1) \).

\eqref{G3D1Y:NL:EF} tells that 
\( \LC{3}{a,b} \) with \( a+b =w, w-1\) are free, and   
\( \LC{3}{a,b} \) with \( a+b \leq w-2\) are determined by 
\( \MC{1}{a,b} \) and  \( \MC{2}{a,b} \) as in  
\eqref{G3D1Y:NL:EF}. Thus, the total number of free parameters is
\( 3\binom{w}{2} + \left( (w+1) + w\right) + \left( 2(w+1) + w-1\right) = 
 3\binom{w}{2} + 5w + 2 = \binom{w}{2} + 2\binom{w+2}{2} + w  \).

Revise on July 20, 2020 by focusing the rank, the number of generators of
linear equations: 
}
\eqref{G3D1Y:NL:E1} $\sim$ \eqref{G3D1Y:NL:EF} yield the generators 
\begin{align*}
& \LC{1}{a,b-1,c+1} (c+1) - \LC{2}{a-1,b,c+1} (c+1)  \\
& \LC{1}{p,q,r+2} 2\tbinom{r+2}{2}  \\
& \LC{2}{p,q,r+2}  2 \tbinom{r+2}{2} \\
& \MC{1}{p-1,q,r+1} (r+1)
- \MC{2}{p,q-1,r+1} (r+1)
+ \LC{3}{p,q,r+2} 2 \tbinom{r+2}{2} 
\end{align*}
The last three are linearly independent and the first type give \(w+1\)
independent equations. Thus the rank is \( 3 \binom{w}{2} +( w +1) \) and the 
kernel dimension is \( 3\binom{w}{2} +  3\binom{w+2}{2} -  3\binom{w}{2}  -
(w+1) = \binom{w}{2} +  2\binom{w+2}{2} +  w\). 

\paragraph{Kernel of 
\(\ds\wtedC{w}{w}\mathop{\to}^{\pdel} \wtedC{w-1}{w} \)}
As a basis of \( \wtedC{w}{w} \), we have 
\(\bU{a,b,c}{w} \) with \( a+b+c= w\) and  
\( \zb{k} \mywedge \zb{4} \mywedge \bU{p,q,r}{w-2} \) with 
\(k=1,2,3\) and \( p+q+r= w-2\),  
take a linear combination
\[ X =  \sum \LC{}{a,b,c} 
\mywedge \bU{a,b,c}{w} 
+ \sum_{k=1}^{3} \sum \MC{k}{p,q,r} z_{k} \mywedge \zb{4} 
\mywedge \bU{p,q,r}{w-2} 
\] by unknown scalars  \(\LC{}{a,b,c},  \MC{k}{p,q,r}  \). 
\begin{align*}
\pdel X =& 
\sum \LC{}{a,b,c} c(c-1) \zb{4} \mywedge \bU{a,b,c-2}{w-2}
+\sum \MC{1}{p,q,r} (-r) \zb{4} \mywedge \bU{p,q+1,r-1}{w-2}
 +\sum \MC{2}{p,q,r} r \zb{4} \mywedge \bU{p+1,q,r-1}{w-2} \;,
\end{align*}
\kmcomment{

we remember \(\MC{3}{p,q,r}\) (\(p+q\leq w-2\)) disappeared here.   
Suppose \( \pdel X=0\). Then we have
\begin{align} \label{G3D1Y:eqn:j1:one}
\sum \LC{}{a,b,c} c(c-1)  \bU{a,b,c-2}{w-2}
+\sum \MC{1}{p,q,r} (-r)  \bU{p,q+1,r-1}{w-2}
 +\sum \MC{2}{p,q,r} r  \bU{p+1,q,r-1}{w-2} = 0 \;
\end{align}
We also remark that  \(\LC{}{a,b,c} \) with \(a+b=w\), 
 \(\LC{}{a,b,c} \) with \(a+b=w-1\), 
 \(\MC{1}{a,b,c} \) with \(a+b=w-2\), and  
 \(\MC{2}{a,b,c} \) with \(a+b=w-2\) disappear in the above equation.   

Shifting the parameters suitably, we reform  \eqref{G3D1Y:eqn:j1:one} as
follows:
\begin{align*}
\sum \LC{}{a,b} (c+2)(c+1)  \bU{a,b}{w-2}
-\sum \MC{1}{a,b-1} (c+1)  \bU{a,b}{w-2}
+\sum \MC{2}{a-1,b} (c+1)  \bU{a,b}{w-2} = 0 \;
\end{align*}
and we get 
\begin{equation}
\LC{}{a,b} = \frac{1}{w-a-b} ( \MC{1}{a,b-1} - \MC{2}{a-1}{b} ) \quad 
( a+b \leq w-2 )\;.
\end{equation}
Thus, the free parameters are just
\( \LC{}{a,b} \) with \( a+b=w, w-1\) and  
\( \MC{i}{a,b} \) with \(i=1,2,3\), \( a+b\leq w-2\). Therefore, the
kernel dimension is \( 2w +1 + 3 \binom{w}{2} = 2\binom{w}{2} + \binom{w+2}
{2} \).

Revise on July 20, 2020 by focusing the rank, the number of generators of
linear equations: 
}
\( \pdel X=0\) is generated by 
\[
\LC{}{p,q,r+2} (r+2)(r+1)  
- \MC{1}{p,q-1,r+1} (r+1)  
 +\MC{2}{p-1,q,r+1} (r+1)  \;, \quad p+q+r = w-2 \;. 
\]
Thus the rank is \(\binom{w}{2}\) and the kernel dimension is 
\( 3\binom{w}{2} + \binom{w+2}{2} - \binom{w}{2} = 
 2\binom{w}{2} + \binom{w+2}{2} \). 

\paragraph{Final table of chain complex }
\begin{thm} \label{D3D1-YinC-src} \ 

\begin{center}
\(
\begin{array}{c|*{5}{c}}
\text{weight}= w & w-1 & w & w+1 & w+2 & w+3\\\hline
\text{SpaceDim} & 
\binom{w}{2} & 
3 \binom{w}{2} + \binom{w+2}{2} & 
3\binom{w}{2} + 3\binom{w+2}{2} & 
\binom{w}{2} + 3\binom{w+2}{2} & 
\binom{w+2}{2} 
\\
\ker\dim & 
\binom{w}{2} & 
 2 \binom{w}{2} + \binom{w+2}{2}  & 
 \binom{w}{2} + 2\binom{w+2}{2} + w  
 &
- \binom{w}{2} + 2 \binom{w+2}{2} & 
w+1 \\\hline
\text{Betti} & 0 & w & 3w+1 & 3w+2  & w+1
\end{array}
\)
\end{center}
\end{thm}


\subsection{\(\dim\frakg=3\), \(\dim [\frakg,\frakg]=1\) and 
\( [\frakg,\frakg] \not \subset Z(\frakg)\) 
}
Consider a Lie algebra \(\frakg\) where 
\([\frakg,\frakg]\) is 1-dimensional and not in of the center of
\(\frakg\). Then we find a basis \(\zb{1}, \zb{2}, \zb{3}\) of \(\frakg\) so
that \(
\Sbt{\zb{1}}{\zb{2}} = -\Sbt{\zb{2}}{\zb{1}} = \zb{2}\) and the
other brackets are zero.
We take \(\ub{1}= \zb{1}\wedge\zb{2},\; \ub{2}= \zb{2}\wedge\zb{3},\;
\ub{3}= \zb{3}\wedge\zb{1}\) as 
a basis of \(\Lambda^{2}\frakg\), and \(V = \zb{4} =  \zb{1} \wedge 
 \zb{2} \wedge  \zb{3} \) as a basis of \(\Lambda^{3}\frakg\). 
Now we have the multiplication (by the Schouten bracket) tables:
\begin{center}
\( \begin{array}[t]{c|*{3}{c}|c|*{3}{c}|}
 & \zb{1} & \zb{2} & \zb{3} & \zb{4} & \ub{1} & \ub{2} & \ub{3} \\\hline 
\zb{1} & 0 & \zb{2} & 0 &\zb{4} & \ub{1} & \ub{2} & 0  \\
\zb{2} & -\zb{2}& 0 & 0 & 0 & 0 & 0 & \ub{2}  \\
\zb{3} & 0 & 0 & 0  & 0 & 0 & 0 & 0 
\end{array}
 \)
\hfil
\( \begin{array}[t]{c|*{3}{c}|c}
 & \ub{1} & \ub{2} & \ub{3} & \zb{4} \\\hline 
\ub{1} & 0 & 0 & \zb{4} & 0 \\
\ub{2} & 0 & 0 & 0 & 0 \\
\ub{3} & \zb{4} & 0 & 0  & 0\\
\end{array}
 \)
\end{center}

Then we have 
\begin{align*}
\pdel W^{1111} &= 2 W^{0111}\;,\quad 
\pdel W^{1110} = \zb{2}\mywedge \zb{3} \;,\quad 
\\
\pdel W^{1101} &= 2 \zb{2} \mywedge \zb{4} \;,\quad 
\pdel W^{1011} =  \zb{3} \mywedge \zb{4}\;,\quad 
\pdel W^{0111} =  0\;,\quad 
\\
\sum_{A} \LC{}{A}\pdel U^{a,b,c} &= \zb{4} \mywedge \sum_{} U^{p,q,r} 
  (p+1)(r+1)  \LC{}{p+1, q, r+1} 
\;. 
\\
\SbtES{\zb{1}}{U^{A}} &=  b U^{a,b-1,c+1} - \beta c U^{a,b+1,c-1}\\ 
\SbtES{\zb{2}}{U^{A}} & = - a U^{a-1,b,c+1} + \alpha c U^{a+1,b,c-1} \\
\SbtES{\zb{3}}{U^{A}} &= \beta a U^{a-1,b+1,c}  - \alpha  b U^{a+1,b-1,c}
\\ \sum_{A} \LC{i}{A} \SbtES{\zb{1}}{U^{A}} 
&= \sum_{A}(a+b)\LC{i} {a,b,c} U^{a,b,c}\;,
\\ \sum_{A} \LC{i}{A} \SbtES{\zb{2}}{U^{A}} 
&= \sum_{A} (c+1)\LC{i}{a,b-1,c+1} U^{a,b,c}\;,
\\  \SbtES{\zb{3}}{U^{A}} 
&= 0 \;. 
\end{align*}

\paragraph{Kernel of \(\ds\wtedC{w+3}{w}\mathop{\to}^{\pdel} \wtedC{w+2}{w} \)}
A basis of \( \wtedC{w+3}{w} \), we have 
\(W^{1,1,1,0 } \mywedge U^{a,b,c} \) with \( a+b+c= w\). 
Take a linear combination
\( X =  \sum_{a+b+c=w} \lambda_{a,b,c} 
W^{1,1,1,0 } \mywedge U^{a,b,c} \) by unknown scalars  \(\lambda_{a,b,c}
\). 
Then  
\begin{align*} \pdel X =&  \sum_{a+b+c= w} 
\lambda_{a,b,c} \left(
(a+b+1) W^{0,1,1,0 } \mywedge U^{a,b,c} 
- c W^{1,0,1,0 } \mywedge U^{a,b+1,c-1} \right.
\\&\qquad \qquad \qquad 
\left. 
- a c 
W^{1,1,1,1 } \mywedge U^{a-1,b,c-1} 
\right)
\\
= &  
 W^{0,1,1,0 }  \sum_{a+b+c= w} \lambda_{a,b,c}
(a+b+1) \mywedge U^{a,b,c} 
-  W^{1,0,1,0 }  \sum_{a+b+c= w} \lambda_{a,b,c} c \mywedge U^{a,b+1,c-1}
\\&\quad 
- W^{1,1,1,1 }   \sum_{a+b+c= w} \lambda_{a,b,c}
 a c   \mywedge U^{a-1,b,c-1} 
\end{align*}
Suppose \(\pdel X = 0\). Then \( \lambda_{a,b,c} (a+b+1) =0 \) 
for \(a+b+c= w\) and all 
\( \lambda_{a,b,c}=0\). Thus, the kernel space is 0-dimensional.  
\kmcomment{
Revise on July 20, 2020 by focusing the rank, the number of generators of
linear equations: 
} Or, 
the linear equations is generated by 
\begin{align*}
& \LC{}{p+1,q,r+1} (p+1)(r+1)\;, \quad ( p+q+r=w-2 ) \\
& \LC{}{a,b,c} (1+a+b)\;,  \quad \LC{}{a,b-1,c+1} (c+1 )\;,  \quad ( a+b+c=w )\;, 
\end{align*}
and the rank is \( \binom{w+2}{2}\) and the kernel dimension is 0.

\paragraph{Kernel of 
\(\ds\wtedC{w+2}{w}\mathop{\to}^{\pdel} \wtedC{w+1}{w} \)}
As a basis of \( \wtedC{w+2}{w} \), we have 
\(
W^{ \eps{1}, \eps{2} ,\eps{3},0} \mywedge \bU{a,b,c}{w}\) with 
\( \eps{1}+ \eps{2} +\eps{3}= 2\) 
and 
\(W^{ 1,1,1 ,1} \mywedge \bU{p,q,r}{w-2}\).   
\kmcomment{
手原稿の \labda^{1} -- [1,1,0] --> \newlabmda^{3}, 
手原稿の \labda^{2} -- [1,0,1] --> \newlabmda^{2}, 
手原稿の \labda^{3} -- [0,1,1] --> \newlabmda^{1}, 
と変更する。1 と 3 の交換のみ。
}

Take a linear combination
\[ X =  \sum_{a+b+c=w} \lambda^{k}_{a,b,c} 
W^{ \eps{1}, \eps{2} ,\eps{3},0} \mywedge \bU{a,b,c}{w}
+  \sum_{p+q+r=w-2} \mu_{p,q,r} 
W^{ 1,1,1,1} \mywedge \bU{p,q,r}{w-2}
\] 
by unknown scalars  \(\lambda^{k}_{a,b,c},  \mu_{p,q,r}  \) 
with 
\( \eps{1}+ \eps{2} +\eps{3}= 2\) and \(\eps{k} = 0\) .
Again, 
\begin{align*}
\pdel X =&  
\sum \lambda^{1}_{a,b,c} ( c W^{0,0,1,0} \bU{a,b+1,c-1}{w}
+ a c  W^{0,1,1,1} \bU{a-1,b,c-1}{w-2})
\\&
+ \sum \lambda^{2}_{a,b,c} ( (a+b) W^{0,0,1,0} \bU{a,b,c}{w}
+ a c  W^{1,0,1,1} \bU{a-1,b,c-1}{w-2})
\\&
+ \sum \lambda^{3}_{a,b,c} ( (a+b+1) W^{0,1,0,0} \bU{a,b,c}{w}
- c W^{1,0,0,0} \bU{a,b+1,c-1}{w}
+ a c W^{1,1,0,1} \bU{a-1,b,c-1}{w-2}
)
\\& 
+ \sum \mu_{p,q,r} ( (p+q+2) W^{0,1,1,1} \bU{p,q,r}{w-2}
- r  W^{1,0,1,1} \bU{p,q+1,r-1}{w-2})
\\ =&  
\sum \lambda^{1}_{a,b,c}  c W^{0,0,1,0} \bU{a,b+1,c-1}{w}
+ \sum \lambda^{2}_{a,b,c} (a+b) W^{0,0,1,0} \bU{a,b,c}{w}
\\&
+ \sum \lambda^{3}_{a,b,c} (a+b+1) W^{0,1,0,0} \bU{a,b,c}{w}
- \sum \lambda^{3}_{a,b,c} (w-a-b) W^{1,0,0,0} \bU{a,b+1}{w}
\\& 
+ \sum \lambda^{1}_{a,b,c}  a(w-a-b) W^{0,1,1,1} \bU{a-1,b,c-1}{w-2}
+ \sum \lambda^{2}_{a,b,c} a(w-a-b) W^{1,0,1,1} \bU{a-1,b,c-1}{w-2}
\\& 
+ \sum \lambda^{3}_{a,b,c} a (w-a-b) W^{1,1,0,1} \bU{a-1,b,c-1}{w-2}
+ \sum \mu_{p,q,r}  (p+q+2) W^{0,1,1,1} \bU{p,q,r}{w-2}
\\& 
- \sum \mu_{p,q,r}  r W^{1,0,1,1} \bU{p,q+1,r-1}{w-2}
\\
=&   W^{0,0,1,0} ( \sum \lambda^{1}_{a,b,c}  c \bU{a,b+1,c-1}{w}
+ \sum \lambda^{2}_{a,b,c} (a+b) \bU{a,b,c}{w}) 
\\&
+W^{0,1,0,0} \sum \lambda^{3}_{a,b,c} (a+b+1) \bU{a,b,c}{w}
-  W^{1,0,0,0} \sum \lambda^{3}_{a,b,c} c \bU{a,b+1,c-1}{w}
\\& 
+  W^{0,1,1,1}(\sum \lambda^{1}_{a,b,c}  a c  \bU{a-1,b,c-1}{w-2}
+ \sum \mu_{p,q,r}  (p+q+2) \bU{p,q,r}{w-2} )
\\& 
+ W^{1,0,1,1} ( \sum \lambda^{2}_{a,b,c} a c  \bU{a-1,b,c-1}{w-2}
-   \sum \mu_{p,q,r}  r  \bU{p,q+1,r-1}{w-2} )
\\& 
+ W^{1,1,0,1} \sum \lambda^{3}_{a,b,c} a c \bU{a-1,b,c-1}{w-2}
\end{align*}

Suppose \( \pdel X=0\). Then we have
\begin{align}
& \sum \lambda^{1}_{a,b,c}  (w-a-b) \bU{a,b+1,c-1}{w}
+ \sum \lambda^{2}_{a,b,c} (a+b) \bU{a,b,c}{w} = 0 
\\&
\sum \lambda^{3}_{a,b,c} (a+b+1) \bU{a,b,c}{w} = 0  \label{G3D1N:lambda3:zero}
\\& 
\sum \lambda^{3}_{a,b,c} (w-a-b) \bU{a,b+1,c-1}{w} = 0 
\\& 
\sum \lambda^{1}_{a,b,c}  a(w-a-b) \bU{a-1,b,c-1}{w-2}
+ \sum \mu_{a,b,c}  (a+b+2) \bU{a,b,c}{w-2} = 0
\\& 
\sum \lambda^{2}_{a,b,c} a(w-a-b) \bU{a-1,b,c-1}{w-2} 
-   \sum \mu_{p,q,r}  r \bU{p,q+1,r-1}{w-2}  = 0
\\& 
 \sum \lambda^{3}_{a,b,c} a c \bU{a-1,b,c-1}{w-2} = 0
\end{align}
\begin{subequations}
The linear equations are generated by 
\begin{align}
& \LC{3}{a,b-1,c+1} (c+1) \label{G3D1N:XX:1} \\
& \LC{3}{a,b,c} (a+b+1) \label{G3D1N:XX:2} \\
& \LC{1}{a,b-1,c+1} (c+1) + \LC{2}{a,b,c} (a+b) \label{G3D1N:XX:3} \\
& \LC{1}{p+1,q,r+1} (p+1)(r+1) + (2+p+q)\MC{ }{p,q,r} \label{G3D1N:XX:4} \\
& \LC{2}{p+1,q,r+1} (p+1) - \MC{ }{p,q-1,r+1} \label{G3D1N:XX:5} \\
& \LC{3}{p+1,q,r+1} (p+1) (r+1) \label{G3D1N:XX:6} 
\end{align}
\end{subequations}
From \eqref{G3D1N:XX:1}, \eqref{G3D1N:XX:2} and \eqref{G3D1N:XX:6}, \(\LC{3}{A}\) are linearly
independent monomials. 

We take linearly independent generators \eqref{G3D1N:XX:3} with leading term
\(\LC{2}{A}\) whose number is \( \binom{w+2}{2} - 1\), where \(a+b=0\) is
the exception.  Finally,   
We take linearly independent generators \eqref{G3D1N:XX:4} with leading term
\(\MC{}{P}\) whose number is \( \binom{w}{2} \).  \eqref{G3D1N:XX:5} are linearly
dependent on \eqref{G3D1N:XX:3} and \eqref{G3D1N:XX:4}. Thus, the rank is 
\( 2 \binom{w+2}{2} - 1  + \binom{w}{2}\) and the kernel dimension is 
\( \binom{w+2}{2} + 1\).

\paragraph{Kernel of 
\(\ds\wtedC{w+1}{w}\mathop{\to}^{\pdel} \wtedC{w}{w} \)}
As a basis of \( \wtedC{w+1}{w} \), we have 
\( \zb{i} \mywedge \bU{a,b,c}{w} \) with \(i=1,2,3\),  and 
\(W^{ \eps{1},\eps{2},\eps{3} ,1} \mywedge \bU{p,q,r}{w-2}\) with   
\( \eps{1}+ \eps{2} +\eps{3}= 2\).  
Take a linear combination 
\[ X = \sum \LC{i}{a,b,c} \zb{i} \mywedge \bU{a,b,c}{w} + 
 \sum_{ \eps{j}=0, \eps{1}+ \eps{2} +\eps{3}= 2 } \MC{j}{p,q,r} 
W^{ \eps{1},\eps{2},\eps{3} ,1} \mywedge \bU{p,q,r}{w-2}\]
for unknown scalars 
\(\LC{i}{a,b,c}\) and \(\MC{j}{p,q,r}\).  
\begin{align*}
\pdel X =& 
\sum \LC{1}{a,b,c} \left(
(a+b)\bU{a,b,c}{w} - a c \zb{1}\mywedge V \mywedge \bU{a-1,b,c-1}{w-2} 
\right)
\\ & 
+ \sum \LC{2}{a,b,c} \left(
c \bU{a,b+1,c-1}{w} - a c \zb{2}\mywedge V \mywedge \bU{a-1,b,c-1}{w-2} 
\right)
\\ & 
- \sum \LC{3}{a,b,c} \left( a c \zb{3}\mywedge V \mywedge \bU{a-1,b,c-1}{w-2} 
\right)
\\ & 
+ \sum \MC{1}{p,q,r} \left( 
r \zb{3}\mywedge V \mywedge \bU{p,q+1,r-1}{w-2} 
\right)
\\ & 
+ \sum \MC{2}{p,q,r} \left( ( p+q-1) \zb{3}\mywedge V \mywedge \bU{p,q,r}{w-2} 
\right)
\\ & 
+ \sum \MC{3}{p,q,r} \left( ( p+q) \zb{2}\mywedge V \mywedge \bU{p,q,r}{w-2} 
  - r \zb{1}\mywedge V \mywedge \bU{p,q+1,r-1}{w-2} 
\right)
\\
=& 
\sum \LC{1}{a,b,c} (a+b)\bU{a,b,c}{w} 
- \sum \LC{1}{a,b,c} a c \zb{1}\mywedge V \mywedge \bU{a-1,b,c-1}{w-2} 
\\ & 
+ \sum \LC{2}{a,b,c} c \bU{a,b+1,c-1}{w} 
- \sum \LC{2}{a,b,c} a c  \zb{2}\mywedge V \mywedge \bU{a-1,b,c-1}{w-2} 
\\ & 
- \sum \LC{3}{a,b,c} a c \zb{3}\mywedge V \mywedge \bU{a-1,b,c-1}{w-2} 
\\ & 
+ \sum \MC{1}{p,q,r} r \zb{3}\mywedge V \mywedge \bU{p,q+1,r-1}{w-2} 
\\ & 
+ \sum \MC{2}{p,q,r} ( p+q-1) \zb{3}\mywedge V \mywedge \bU{p,q,r}{w-2} 
\\ & 
+ \sum \MC{3}{p,q,r} ( p+q) \zb{2}\mywedge V \mywedge \bU{p,q,r}{w-2} 
- \sum \MC{3}{p,q,r} r \zb{1}\mywedge V \mywedge \bU{p,q+1,r-1}{w-2} 
\\ =& 
\sum \LC{1}{a,b,c} (a+b)\bU{a,b,c}{w} 
+ \sum \LC{2}{a,b,c} c \bU{a,b+1,c-1}{w} 
\\ & - \zb{1}\mywedge V \mywedge (
 \sum \LC{1}{a,b,c} a c   \bU{a-1,b,c-1}{w-2} 
+ \sum \MC{3}{p,q,r} r \bU{p,q+1,r-1}{w-2} 
)
\\& + \zb{2}\mywedge V \mywedge(
 \sum \MC{3}{p,q,r} ( p+q)  \bU{p,q,r}{w-2} 
- \sum \LC{2}{a,b,c} a c \bU{a-1,b,c-1}{w-2} )
\\ & + \zb{3}\mywedge V \mywedge (
- \sum \LC{3}{a,b,c} a c  \bU{a-1,b,c-1}{w-2} 
 + \sum \MC{1}{p,q,r}  r )  \bU{p,q+1,r-1}{w-2} 
+ \sum \MC{2}{p,q,r} ( p+q-1) \bU{p,q,r}{w-2} )
\end{align*}
Thus, we have linear equations 
\begin{align}
& \sum \LC{1}{a,b,c} (a+b)\bU{a,b,c}{w} 
+ \sum \LC{2}{a,b,c} c \bU{a,b+1,c-1}{w}  = 0 \tag{E1} \label{G3D1N:L:E1}
\\ & \sum \LC{1}{a,b,c} a c   \bU{a-1,b,c-1}{w-2} 
+ \sum \MC{3}{p,q,r} r \bU{p,q+1,r-1}{w-2} 
= 0 \tag{E2} \label{G3D1N:L:E2}
\\& \sum \MC{3}{p,q,r} ( p+q)  \bU{p,q,r}{w-2} 
- \sum \LC{2}{a,b,c} a c \bU{a-1,b,c-1}{w-2} = 0 \tag{E3} \label{G3D1N:L:E3}
\\ & \sum \LC{3}{a,b,c} a c   \bU{a-1,b,c-1}{w-2} 
- \sum \MC{1}{p,q,r}  r   \bU{p,q+1,r-1}{w-2} 
- \sum \MC{2}{p,q,r} ( p+q-1) \bU{p,q,r}{w-2} = 0 \tag{EF} \label{G3D1N:L:EF}
\end{align}
\begin{subequations}
From \eqref{G3D1N:L:E1} $\sim$ \eqref{G3D1N:L:EF}, we have generators 
\begin{align}
& \LC{1}{a,b,c} (a+b) + \LC{2}{a,b-1,c+1} (c+1) \label{G3D1N:Y:1}
\\& \LC{1}{p+1,q,r+1} (p+1)(r+1) +  \MC{3}{p,q-1,r+1} (r+1)  \label{G3D1N:Y:2} 
\\&  \LC{2}{p+1,q,r+1} (p+1)(r+1) - \MC{3}{p,q,r} (2+p+q ) \label{G3D1N:Y:3}
\\&  \LC{3}{p+1,q,r+1} (p+1)(r+1) 
- \MC{1}{p,q-1,r+1} (r+1) - \MC{2}{p,q,r} ( p+q+1)  \label{G3D1N:Y:4}
\end{align}
\end{subequations}
We take linearly independent generators 
\eqref{G3D1N:Y:1} with leading term \( \LC{1}{a,b,c}(a+b) \), 
\eqref{G3D1N:Y:3} with leading term \( \MC{3}{p,q,r} \),  and 
\eqref{G3D1N:Y:4} with leading term \( \MC{2}{p,q,r} \).  Then \eqref{G3D1N:Y:2} are
linearly dependent on those of 3 types. Thus, the rank is 
\( \binom{w+2}{2} - 1 + 2 \binom{w}{2}\) and the kernel dimension is 
\( 1+\binom{w}{2} + 2 \binom{w+2}{2} \). 

\paragraph{Kernel of 
\(\ds\wtedC{w}{w}\mathop{\to}^{\pdel} \wtedC{w-1}{w} \)}
As a basis of \( \wtedC{w}{w} \), we have 
\(U^{a,b,c} \) with \( a+b+c= w\) and  
\( z_{k} \mywedge V \mywedge U^{p,q,r} \) with 
\(k=1,2,3\) and \( p+q+r= w-2\). 
\kmcomment{
To avoid confusion, we use the
notation \(\ds \bU{a,b,c}{\ell}\) which  means 
\(\ds U^{a,b,c}\) with the restriction \(a+b+c= \ell\). The third
component \(c \) is uniquely determined by \(c =\ell-a-b\), we may
denote \(\ds \bU{a,b,c}{\ell}\) by 
\(\ds \bU{a,b}{\ell}\) with \( a+b \leqq \ell\).  
Using the notations explained in Remark \ref{G3D1N:remark:odd:nota}, 
}
A basis of \( \wtedC{w}{w} \) is 
\(\bU{a,b,c}{w} \) and \( z_{k} \mywedge V \mywedge \bU{p,q,r}{w-2} \). 
Consider a linear combination
\[ X =  \sum_{a+b+c=w} \LC{}{a,b,c} \bU{a,b,c}{w} 
+ \sum_{k=1}^{3} \sum_{p+q+r=w-2}  \MC{k}{p,q,r} z_{k} 
\mywedge V 
\mywedge \bU{p,q,r}{w-2} 
\] by unknown scalars  \(\LC{}{a,b,c},  \MC{k}{p,q,r}  \). 
Then 
\begin{align*}
\pdel X =& 
\sum \LC{}{a,b,c} a c V \mywedge \bU{a-1,b,c-1}{w-2}
+\sum \MC{1}{p,q,r} (p+q+1) V \mywedge \bU{p,q,r}{w-2}
\\& 
+\sum \MC{2}{p,q,r} r V \mywedge \bU{p,q+1,r-1}{w-2} \;.
\end{align*}
\kmcomment{
Suppose \( \pdel X=0\). Looking at \( V\mywedge \bU{a,b}{w-2}\) when
\(b=0\) or when positive \(b\), we have  
 \begin{align*}
 & \lambda_{1+a,0}(1+a)(w-1-a) + \mu_{1,a,0} (a+1) = 0 \\
 & \lambda_{1+a,b}(1+a)(w-1-a-b) + 
 \mu_{1,a,b}(a+b+1) + 
 \mu_{2,a,b-1} (w-3-a-b) = 0 \quad ( b>0 )
 \\\noalign{and so}
 & \mu_{1,a,0}  = - (w-1-a)  \lambda_{1+a,0} \\ 
 & 
 \mu_{1,a,b}  = - 
\frac{ (1+a)(w-1-a-b) }{a+b+1}  \lambda_{1+a,b}
 - \frac{ (w-3-a-b) }{(a+b+1) }
 \mu_{2,a,b-1} \quad ( b>0 )\;.
\end{align*}
Thus, the free parameters are \( \lambda_{a,b}\) with \( a+b+c=w\), and 
\( \mu_{2,a,b}\) and 
\( \mu_{3,a,b}\) with \( a+b+c=w-2\).  The kernel dimension is \(
2\binom{w}{2} + 
\binom{w+2}{2}\).   

Revise on July 20, 2020 by focusing the rank, the number of generators of
linear equations: 
}
Since \( \pdel X = V \mywedge \sum \bU{p,q,r}{w-2} ( \LC{}{p+1,q,r+1} (p+1)
(r+1) + (1+p+q) \MC{1}{p,q,r} +  (1+r) \MC{2}{p,q-1, r+1}  ) \), the rank is
\( \binom{w}{2} \)  because the leading terms are \( \MC{1}{p,q,r}\), thus
the kernel dimension is \( ( 3\binom{w}{2} + \binom{w+2}{2} ) - \binom{w}{2}
= 2\binom{w}{2} + \binom{w+2}{2}  \).

\paragraph{Final table of chain complex }
\begin{thm} \label{thm:G3D1-NinC} \ 

\begin{center}
\(
\begin{array}{c|*{5}{c}}
\text{weight}= w & w-1 & w & w+1 & w+2 & w+3\\\hline
\text{SpaceDim} & 
\binom{w}{2} & 
3 \binom{w}{2} + \binom{w+2}{2} & 
3\binom{w}{2} + 3\binom{w+2}{2} & 
\binom{w}{2} + 3\binom{w+2}{2} & 
\binom{w+2}{2} 
\\
\ker\dim & 
\binom{w}{2} & 
2 \binom{w}{2} + \binom{w+2}{2} & 
1+ \binom{w}{2} + 2\binom{w+2}{2} & 
1+\binom{w+2}{2} & 
0 \\\hline
\text{Betti} & 0 & 1 & 2 & 1 & 0
\end{array}
\)
\end{center}
\end{thm}

\kmcomment{
\begin{titlepage} 
\vspace*{20mm}
\begin{center}
{\Large Super homologies associated with low dimensional Lie algebras}

\vspace{10mm}
{Kentaro Mikami ({Akita University}) 
 \quad Tadayoshi Mizutani ({Saitama University})}

\vspace{10mm}
{\today} 
\end{center}
\end{titlepage}
}

\subsection{\(\dim[\frakg,\frakg]=2\)}
\renewcommand{\arraystretch}{1.0}
Take a 3-dimensional Lie algebra whose derived subalgebra is 2-dimensional.
Then we have a basis 
satisfying 
\[ 
\Sbt{\zb{1}}{\zb{2}} =0,\  \Sbt{\zb{1}}{\zb{3}} = \zb{1},\ 
\Sbt{\zb{2}}{\zb{3}} = \alpha \zb{2} \quad( \alpha \ne 0 ) \; .
\] 
Now we have the multiplication (by the Schouten bracket) tables with
notations  
\( \ub{1} = \zb{1} \wedge \zb{2}\;,\;  \ub{2} = \zb{1} \wedge \zb{3}\;,\;  
\ub{3} = \zb{2} \wedge \zb{3}\;,\; 
\zb{4} = V = \zb{1} \wedge \zb{2} \wedge \zb{3}\).
\begin{center}
\( \begin{array}[t]{c|*{3}{c}|c|*{3}{c}|}
 & \zb{1} & \zb{2} & \zb{3} & \zb{4} & \ub{1} & \ub{2} & \ub{3}  
 \\\hline 
\zb{1} & 0 &  0 & \zb{1} &0  & 0 & 0 & -\ub{1}   \\
\zb{2} &  0 & 0 & \alpha \zb{2}& 0 &  0 & \alpha \ub{1} & 0  \\
\zb{3} & - \zb{1} & -\alpha \zb{2} & 0 & -(1+\alpha) V & -(1+\alpha) \ub{1}  & -\ub{2}  &
- \alpha\ub{3} 
\end{array}
 \)
\hfil
\( \begin{array}[t]{c|*{3}{c}|c}
   & \ub{1} & \ub{2} & \ub{3} & \zb{4} 
\\\hline 
\ub{1} & 0 & 0 & 0 & 0 \\
\ub{2} & 0 & 0 &(1-\alpha)\zb{4}  & 0 \\
\ub{3} & 0 & (1-\alpha)\zb{4} &  0  & 0\\
\end{array}
 \)
\end{center}

Then we have 
\begin{align*}
\pdel W^{1111} &= -2 (1+\alpha) \zb{1}\mywedge  \zb{2}\mywedge\zb{4}
\;,\quad 
\pdel W^{1110} =  - (1+\alpha)  \zb{1}\mywedge\zb{2}\;,\quad 
\\
\pdel W^{1101} &= 0 \;,\quad 
\pdel W^{1011} =  (2+\alpha) \zb{1} \mywedge \zb{4}\;,\quad 
\pdel W^{0111} = (2 \alpha +1) \zb{2} \mywedge \zb{4}\;,\quad 
\\
\pdel U^{a,b,c} &= \zb{4} \mywedge (
  b c  U^{a, b-1, c-1} 
)
\;, \quad  
\sum_{A} \LC{i}{A}\pdel U^{a,b,c} = \zb{4} \mywedge \sum_{i} U^{p,q,r} (
  (q+1)(r+1)  \LC{i}{p, q+1, r+1} 
)
\;. 
\\
\SbtES{\zb{1}}{U^{A}} &=  -  c U^{a+1,b,c-1}\;,\quad  
\sum_{A} \LC{i}{A} \SbtES{\zb{1}}{U^{A}} 
= \sum_{A}(-(c+1)\LC{i}{a-1,b,c+1}) U^{A}\;,
\\ 
\SbtES{\zb{2}}{U^{A}} & = \alpha b U^{a+1,b-1,c} \;,\quad
\sum_{A} \LC{i}{A} \SbtES{\zb{2}}{U^{A}} 
= \sum_{A}(\alpha (b+1)\LC{i}{a-1,b+1,c}) U^{A}\;,
\\ 
\SbtES{\zb{3}}{U^{A}} &= ( -(1+\alpha) a -b -\alpha c ) U^{a,b,c} 
\; . 
\end{align*}

\begin{remark}\label{G3D2:remark:parameter}
When the weight is 0, the usual Lie algebra homology groups are
obtained as below: 
\begin{center}
\(
\begin{array}{c|*{4}{c}}
& \Lambda^{0} \frakg & \Lambda^{1} \frakg & \Lambda^{2} \frakg 
        & \Lambda^{3} \frakg \\\hline
\text{SpaceDim} & 1 & 3 & 3 & 1 \\
\ker \dim  & 1 & 3 & 1 & \kappa \\\hline
\text{Betti} & 1 & 1 & \kappa & \kappa \\\hline
\end{array}
\)
\end{center}
Since \( \pdel ( \zb{1} \wedge  \zb{2} \wedge  \zb{3}) = - (1+\alpha)
\zb{1} \wedge \zb{2} \), we have 
\( \kappa = \begin{cases} 1 & \text{if\quad} \alpha = -1 \\
0 & \text{if\quad} \alpha \ne  -1 
\end{cases}
\). Encountering with this phenomena, we are interested in 
how the super homology groups depend on non-zero \(\alpha\), which
is a parameter of  the Lie algebra structures. 
\end{remark}

\paragraph{The space of cycles in \(\wtedC{w+3}{w}\):} 
We study the space of cycles in \(\wtedC{w+3}{w}\). 
Take \( W^{1,1,1,0} \mywedge U^{a,b,c}\)
with \(a+b+c=w\) as a  basis of  the chain space.   
Consider a general chain 
\( \sum_{a+b+c=w} \LC{}{a,b,c}  W^{1,1,1,0} \mywedge U^{a,b,c}\) with
unknown scalars \(\LC{}{a,b,c}\). 
\begin{align*} 
& \pdel ( \sum_{a+b+c=w} \LC{}{a,b,c}  W^{1,1,1,0} \mywedge U^{a,b,c})
\\ =& 
 W^{0,1,1,0} \mywedge \sum_{a+b+c=w} \LC{}{a,b,c}  c U^{a+1,b,c-1}
-  W^{1,0,1,0} \mywedge \sum_{a+b+c=w} \LC{}{a,b,c} \alpha  b
U^{a+1,b-1,c}
\\& 
- W^{1,1,0,0} \mywedge \sum_{a+b+c=w} \LC{}{a,b,c} ( (1+\alpha)( a+1)  + b + \alpha c ) 
U^{a,b,c}
\\& 
- W^{1,1,1,1} \mywedge \sum_{a+b+c=w} \LC{}{a,b,c}  bc(1-\alpha) 
U^{a,b-1,c-1}
\end{align*} 
Now assume that the above is zero.  Then we have
\begin{subequations}
\begin{align}
& 
 \sum_{a+b+c=w} \LC{}{a,b,c}  c U^{a+1,b,c-1} = 0 \label{G3D2:z:e1}
\\& 
 \sum_{a+b+c=w} \LC{}{a,b,c}  b U^{a+1,b-1,c} = 0\label{G3D2:z:e2}
\\& 
 \sum_{a+b+c=w} \LC{}{a,b,c} 
 ( (1+\alpha)(1+a)  + b + \alpha c ) U^{a,b,c}= 0 \label{G3D2:z:e3}
\\& 
 \sum_{a+b+c=w} \LC{}{a,b,c}  bc(1-\alpha) 
U^{a,b-1,c-1} = 0\label{G3D2:z:e4}
\end{align}
\kmcomment{
\eqref{G3D2:z:e1} and 
\eqref{G3D2:z:e2} imply \( \LC{}{a,b,c} =0\) for \(c>0\) or  \(b>0\). When
\(b=c=0\),  \(a=w\) and \eqref{G3D2:z:e3} shows that
\begin{equation} \label{G3D2:z:e5}
 \LC{}{w,0,0}   (1+\alpha) (1+w)  U^{w,0,0}= 0
\end{equation}

If \(w=0\), then we have 
\(
 \LC{}{0,0,0} ( 1+\alpha ) U^{0,0,0}= 0
 \) and this explains the case of Remark 
\ref{G3D2:remark:parameter}. 

Even when \(w>0\), \eqref{G3D2:z:e5} says 
the kernel dimension of \(\wtedC{w+3}{w}\) is 1 if \( \alpha = -1
\), otherwise that is 0.      

\textcolor{red}{Revised} on July 21, 2020 by focusing the rank, the number of generators of linear equations: 
}
\eqref{G3D2:z:e1} $\sim$ 
\eqref{G3D2:z:e4} yield generators of linear equation system as follows:
\begin{align}
& \LC{}{a-1,b,c+1}  (c+1) \label{G3D2:z2:e1}
\\& \LC{}{a-1,b+1,c}  (b+1)\label{G3D2:z2:e2}
\\& \LC{}{a,b,c} ( (1+\alpha)(1+a)  + b + \alpha c ) \label{G3D2:z2:e3}
\\&  \LC{}{a,b+1,c+1}  (b+1)(c+1)(1-\alpha) \label{G3D2:z2:e4}
\end{align}
If \(\alpha=1\), then we see directly \(\LC{}{A} \) are the linearly
independent generators by \eqref{G3D2:z2:e3}, and so the rank is \(\binom{w+2}{2}
\). 

If \(\alpha \ne 1\), 
\eqref{G3D2:z2:e4} provide \(\LC{}{a,b,c}\) with \(b c>0\),   
\eqref{G3D2:z2:e1} provide \(\LC{}{a,0,c}\) with \( c>0\),   
\eqref{G3D2:z2:e2} provide \(\LC{}{a,b,0}\) with \( b>0\),   and 
\eqref{G3D2:z2:e3} provide \(\LC{}{a,0,0}(1+\alpha)(1+a)\).   
Thus, the rank is \(\binom{w+2}{2}-1 \)  if \( \alpha+1=0\) or 
 the rank is \(\binom{w+2}{2} \)  if \( \alpha+1 \ne 0\).  
\end{subequations}
\bigskip

\paragraph{The space of cycles in \(\wtedC{w+2}{w}\):} 
We study the space of  cycles in \(\wtedC{w+2}{w}\). 
Take 
\( W[\eps{i}=0,\eps{4}=0] \mywedge \bU{a,b,c}{w}\)
for \(i=1,2,3\) and \( W^{1,1,1,1} \mywedge \bU{p,q,r}{w-2}\) as a basis of    
the chain space, and   consider a general chain 
\( \sum \LC{i}{a,b,c}  W[{\eps{i}=0,\eps{4}=0}] \mywedge \bU{a,b,c}{w}
+ \sum \MC{} {p,q,r}  W^{1,1,1,1} \mywedge \bU{p,q,r}{w-2}
\) with unknown scalars \(\LC{i}{a,b,c}, \MC{}{p,q,r}\).   
\begin{align*}
&\ \pdel (
 \sum \LC{i}{a,b,c}  W[{\eps{i}=0,\eps{4}=0}] \mywedge \bU{a,b,c}{w}
+ \sum \MC{} {p,q,r}  W^{1,1,1,1} \mywedge \bU{p,q,r}{w-2}) 
\\ =&\    
+ W^{1,0,0,0} \mywedge \sum ( \LC{2}{a,b,c}( (1+\alpha)a + b+ \alpha c +1)
\bU{a,b,c}{w} - \LC{3}{a,b,c} \alpha b  \bU{a+1,b-1,c}{w}) 
\\& 
+ W^{0,1,0,0} \mywedge\sum ( \LC{1}{a,b,c}( (1+\alpha)a + b+ \alpha c
+ \alpha )
\bU{a,b,c}{w} - \LC{3}{a,b,c} c \bU{a+1,b,c-1}{w}) 
\\& 
+ W^{0,0,1,0} \mywedge\sum ( \alpha \LC{1}{a,b,c} b \bU{a+1,b-1,c}{w} 
- \LC{2}{a,b,c} c \bU{a+1,b,c-1}{w}) 
\\& 
+ W^{0,1,1,1} \mywedge \left( - \sum  \MC{ }{p,q,r} r \bU{p+1,q,r-1}{w-2} 
+ \sum \LC{1}{a,b,c} (1-\alpha) b c \bU{a,b-1,c-1}{w-2}    
\right)
\\& 
+ W^{1,0,1,1} \mywedge\left(  \sum  \MC{ }{p,q,r}(- \alpha q) \bU{p+1,q-1,r}{w-2} 
+ \sum \LC{2}{a,b,c} (1-\alpha) b c \bU{a,b-1,c-1}{w-2}    
\right)
\\& 
+ W^{1,1,0,1} \mywedge\left(\sum \MC{ }{p,q,r}(- (1+\alpha)(p+2)+ q+ \alpha r )
\bU{p,q,r}{w-2} + \sum  \LC{3}{a,b,c} (1-\alpha) b c \bU{a,b-1,c-1}{w-2}    
\right)
\end{align*}

We pick up the coefficient of \(W^{\varE} \mywedge U^{a,b,c}\):   
\begin{subequations}
\begin{align}
 W^{1,0,0,0} & : \quad  \LC{2}{a,b,c}(  (1+\alpha)a + b+ \alpha c + 1
) - \alpha \LC{3}{a-1,b+1,c} (b+1)   \label{G3D2:f:0:a}
\\ 
 W^{0,1,0,0} & : \quad  \LC{1}{a,b,c}(
 (1+\alpha)a + b+ \alpha c  + \alpha  )
 - \LC{3}{a-1,b,c+1} (c+1)  \label{G3D2:f:0:b}
\\ 
 W^{0,0,1,0}& : \quad   \alpha  \LC{1}{a-1,b+1,c}  (b+1) 
- \LC{2}{a-1,b,c+1} (c+1)  \label{G3D2:f:0:c}
\\ 
 W^{0,1,1,1} & : \quad -  \MC{ }{p-1,q,r+1} (r+1) 
 + (1-\alpha) \LC{1}{p,q+1,r+1} (q+1) (r+1) 
 \label{G3D2:f:0:d}
\\ 
 W^{1,0,1,1} & : \quad  
 - \alpha  \MC{ }{p-1,q+1,r} (q+1) 
 + (1-\alpha) \LC{2}{p,q+1,r+1} (q+1) (r+1) 
\label{G3D2:f:0:e}
\\ 
 W^{1,1,0,1}& :  \quad  
- \MC{ }{p,q,r}( (1+\alpha)(p+2) + q+ \alpha r  )
 + (1-\alpha) \LC{3}{p,q+1,r+1} (q+1) (r+1) 
 \label{G3D2:f:0:f}
\end{align}
\end{subequations}
If \( \alpha=1\), then \eqref{G3D2:f:0:d} $\sim$ \eqref{G3D2:f:0:f} yield independent
generators \( \MC{}{p,q,r}\) with \( p+q+r=w-2\). \eqref{G3D2:f:0:a} and 
\eqref{G3D2:f:0:b} are independent generators which express \eqref{G3D2:f:0:c}. Thus, the rank
is \( \binom{w}{2} + 2 \binom{w+2}{2}\). 

If \( \alpha\ne 1\),  \eqref{G3D2:f:0:d} $\sim$ \eqref{G3D2:f:0:f} yield independent
generators. Be careful that \eqref{G3D2:f:0:f} with \(q=r=0\) is 
\( - \MC{}{p,0,0}(1+\alpha)(p+2) + (1-\alpha) \LC{3}{p,1,1}\). 
We need more generators from 
\eqref{G3D2:f:0:a} $\sim$ \eqref{G3D2:f:0:c} with the condition \( bc=0\). 
\eqref{G3D2:f:0:a} with \(b=0, c=1\) is 
\( -\LC{2}{w-1,0,1}(1+\alpha) w -\alpha \LC{3}{w-2,1,1}\). 
\eqref{G3D2:f:0:b} with \(b=1, c=0\) is 
\( \LC{1}{w-1,1,0}(1+\alpha) w - \LC{3}{w-2,1,1}\). 

If \( \alpha + 1 \ne 0\), then \eqref{G3D2:f:0:c} are expressed by 
\eqref{G3D2:f:0:a} and \eqref{G3D2:f:0:b}, the rank is \( 3\binom{w}{2} + 2(2w+1) =
\binom{w}{2} + 2 \binom{w+2}{2}\).  

If \(\alpha + 1 = 0\), as described above, we have 3 cases which mean the
same generator \(\LC{3}{w-2,1,1}\). Depending to those cases, we need
relation \( \alpha \LC{1}{w-1,1,0} - \LC{2}{w-1,0,1}\) from \eqref{G3D2:f:0:c}. 
Thus, the rank is \( 3 \binom{w}{2} + 2 (w+1) - 2 + 1 = 
\binom{w}{2} + 2 \binom{w+2}{2} -1 \).

\paragraph{The space of cycles in \(\wtedC{w+1}{w}\):} 
When \(w=0\), then \(\wtedC{w+1}{w} = \wtedC{1}{0} = \frakg\) and the
boundary operator is trivial, and so the kernel dimension is \(\dim
\frakg\). 
Thus, we study the space of cycles in  \(\wtedC{w+1}{w}\) for \(w > 0\).  
Take \( \zb{i} \mywedge \bU{a,b,c}{w}\) and  \( \zb{i} \mywedge  \zb{j}
\mywedge \zb{4} \mywedge \bU{a,b,c}{w-2}\) (\(1\leqq i<j\leqq3\)) as a
basis of the chain space, and      consider a general chain \(L+M\) where \(
L = \sum \LC{i}{a,b,c}  \zb{i} \mywedge \bU{a,b,c}{w}\) and \( M = \sum
\MC{i} {p,q,r} W[\eps{i}=0, \eps{4}=1 ]\mywedge \bU{p,q,r}{w-2} \) with
unknown scalars \(\LC{i}{a,b,c}, \MC{i}{p,q,r}\).   
\begin{align*}
\pdel L =& 
- \zb{1}\mywedge \sum \LC{1}{a,b,c}\pdel\bU{a,b,c}{w}- \sum \LC{1}{a,b,c}
c\bU{a+1,b,c-1}{w} \\&
- \zb{2}\mywedge \sum \LC{2}{a,b,c}\pdel\bU{a,b,c}{w}+ \sum \LC{2}{a,b,c} \alpha
b\bU{a+1,b-1,c}{w} \\&
- \zb{3}\mywedge \sum \LC{3}{a,b,c}\pdel\bU{a,b,c}{w}- \sum \LC{3}{a,b,c} 
( (1 +\alpha) a + b + \alpha c ) \bU{a,b,c}{w} 
\end{align*} 
Using the fact \(
\SbtES{\zb{i}\mywedge \zb{j}\mywedge \zb{4} }{U^{a,b,c}} = 
\zb{j} \mywedge \zb{4} \mywedge \SbtES{\zb{i}}{ U^{a,b,c}}
- \zb{i} \mywedge \zb{4} \mywedge \SbtES{\zb{j}}{ U^{a,b,c}}
\)  we have
\begin{align*}
\pdel M =& 
+ W^{0,1,0,1} \mywedge  \sum \MC{1}{P} 
( 1+2 \alpha + (1+\alpha) p + q + \alpha r) U^{P} 
+ W^{0,0,1,1} \mywedge  \sum \MC{1}{P} \alpha q U^{p+1,q-1,r}
\\&
+ W^{1,0,0,1} \mywedge  \sum \MC{2}{P} 
( 2+\alpha + (1+\alpha) p + q + \alpha r) U^{P} 
+ W^{0,0,1,1} \mywedge  \sum \MC{2}{P} ( -r ) U^{p+1,q,r-1}
\\&
+ W^{1,0,0,1} \mywedge  \sum \MC{3}{P} 
( - \alpha r U^{p,q+1,r-1} ) 
+ W^{0,1,0,1} \mywedge  \sum \MC{3}{P} ( -r ) U^{p+1,q,r-1}
\end{align*}
Assume that \(\pdel(L +  M) = 0\). 
Knowing \( \pdel U^{a,b,c} = \zb{4} \mywedge (1-\alpha) b c U^{a,b-1,c-1}\), 
\begin{subequations}
\begin{align*} 
0 & = \sum \LC{1}{a,b,c} (1-\alpha) b c \bU{a,b-1,c-1}{w} 
- \sum \MC{2}{P} ( 2+\alpha + (1+\alpha) p + q + \alpha r) U^{P} 
\\& \qquad 
+  \sum \MC{3}{P} (  \alpha r ) U^{p,q+1,r-1} 
\notag
\\ 0 &  =   
\sum \LC{2}{a,b,c} (1-\alpha) b c \bU{a,b-1,c-1}{w} 
- \sum \MC{1}{P} 
( 1+2\alpha + (1+\alpha) p + q + \alpha r) U^{P} 
\\& \qquad 
-   \sum \MC{3}{P} ( -r ) U^{p+1,q,r-1}
\notag
\\0 & =  
\sum \LC{3}{a,b,c} (1-\alpha) b c \bU{a,b-1,c-1}{w} 
-   \sum \MC{1}{P} \alpha q U^{p+1,q-1,r}
-  \sum \MC{2}{P} ( -r ) U^{p+1,q,r-1}
\\0 & =   - 
\sum \LC{1}{a,b,c} c\bU{a+1,b,c-1}{w} 
+ \sum \LC{2}{a,b,c} \alpha b\bU{a+1,b-1,c}{w}  
\\& \qquad 
- \sum \LC{3}{a,b,c} ( (1 +\alpha) a + b + \alpha c ) \bU{a,b,c}{w} 
\notag
\end{align*}
Taking coefficient of \( \bU{a,b,c}{}\), we have
\begin{align} \label{G3D2:G:0:a}
 &   (1-\alpha)\LC{1}{p,q+1,r+1} (q+1)( r+1) 
- \MC{2}{p,q,r} ( 2+\alpha + (1+\alpha) p + q + \alpha r) 
\\& \qquad 
+ \alpha  \MC{3}{p,q-1,r+1} (r+1) 
\notag
\\  &   
 (1-\alpha) \LC{2}{p,q+1,r+1}(q+1)(r+1) 
-  \MC{1}{p,q,r} 
( 1+2\alpha + (1+\alpha) p + q + \alpha r) 
\label{G3D2:G:0:b}
\\& \qquad 
+   \MC{3}{p-1,q,r+1} ( r+1 )
\notag
\\ &   \label{G3D2:G:0:c}
(1-\alpha) \LC{3}{p,q+1,r+1} (q+1)(r+1) 
- \alpha  \MC{1}{p-1,q+1,r} (q+1)
+ \MC{2}{p-1,q,r+1} ( r+1 )
\\  &   -
 \LC{1}{a-1,b,c+1} (c+1)
+  \alpha \LC{2}{a-1,b+1,c} (b+1)   \label{G3D2:G:0:d}
-  \LC{3}{a,b,c} ( (1 +\alpha) a + b + \alpha c ) 
\end{align}
\end{subequations}
Assume \( 1-\alpha = 0\). 
\begin{subequations}
We have
\begin{align}
 & \label{G3D2:G:k:a}
- \MC{2}{p,q,r} ( 1 + w + p ) 
 +   \MC{3}{p,q-1,r+1}  (r+1) 
\\  &    \label{G3D2:G:k:b}
-  \MC{1}{p,q,r} ( 1 + w + p ) 
+    \MC{3}{p-1,q,r+1} (r+1)  
\\ &    \label{G3D2:G:k:c}
-   \MC{1}{p-1,q+1,r}  (q+1)
+  \MC{2}{p-1,q,r+1} (r+1 )
%
\\
 &   
 - \LC{1}{a-1,b,c+1} (c+1)
+  \LC{2}{a-1,b+1,c} (b+1) 
-  \LC{3}{a,b,c} ( w + a ) \label{G3D2:G:k:d}
\end{align}
\end{subequations}
When \( \alpha = 1\), from \eqref{G3D2:G:0:a} and \eqref{G3D2:G:0:b}, we have
\(2\binom{w}{2}\) generators and \eqref{G3D2:G:0:c}  are controlled by  
\eqref{G3D2:G:0:a} and \eqref{G3D2:G:0:b}. 
\eqref{G3D2:G:0:d} yield \(\binom{w+2}{2}\) generators. Thus, the rank is \(
2\binom{w}{2} + \binom{w+2}{2}\).  

If 
\( \alpha \ne  1\), from \eqref{G3D2:G:0:a} $\sim$ \eqref{G3D2:G:0:c}, we have
\(3 (  \binom{w+2}{2}- (2w+1)) \) generators,  and \eqref{G3D2:G:0:d}  yield 
generators depending on \((a,b,c)\) with \( b c =0\). Thus, the rank is 
\( 3( \binom{w+2}{2} - (2w+1) ) + (2w+1) = 
3\binom{w+2}{2} - 2(2w+1) = \binom{w+2}{2} + 2\binom{w}{2} \).

\paragraph{The space of cycles in \(\wtedC{w}{w}\):} 
When \(w=0\), then \(\wtedC{w}{w} = \Lambda ^{0}\frakg\) and the
boundary operator is trivial, and so the kernel dimension is 1.  Thus,
we study the space of cycles in  \(\wtedC{w}{w}\)) when \(w > 0\).  
Take 
\( \bU{a,b,c}{w}\) and  \( \zb{i} \mywedge \zb{4} \mywedge
\bU{p,q,r}{w-2}\) (\(i=1,2,3\)) as a basis of the chain
space, and  
consider a general chain \(L+M\) where \( L = \sum \LC{ }{a,b,c}
\bU{a,b,c}{w}\) and \( M = \sum \MC{i} {p,q,r}  \zb{i} \mywedge \zb{4}
\mywedge \bU{p,q,r}{w-2} \) with unknown scalars \(\LC{ }{a,b,c},
\MC{i}{p,q,r}\).   
\begin{align*}
\pdel L =& 
(1-\alpha) \zb{4} \mywedge \sum  \LC{ }{a,b,c}  b c \bU{a,b-1,c-1}{w-2} 
= (1-\alpha) \zb{4} \mywedge \sum  \LC{ }{a,b+1,c+1} (b+1)(c+1)
\bU{a,b,c}{w-2} 
\end{align*}
Using the facts \( \zb{4}\mywedge \pdel U^{A}=0\) and  \(
\SbtES{\zb{i}\mywedge \zb{4}}{U^{A}} = \zb{4} \mywedge 
\SbtES{\zb{i}}{ U^{A}}\),  we have

\begin{align*}
\pdel M =& 
-\zb{4}\mywedge   \sum \MC{1}{p,q,r} c\bU{p+1,q,r-1}{w-2} 
+\zb{4} \mywedge \sum \MC{2}{p,q,r} \alpha q\bU{p+1,q-1,r}{w-2} \\&
 - \zb{4}\mywedge  \sum \MC{3}{p,q,r} (1+\alpha) \bU{p,q,r}{w-2}
-\zb{4}\mywedge \sum \MC{3}{p,q,r} ( (1 +\alpha) p + q + \alpha r ) \bU{p,q,r}{w-2} 
\\
=& 
\zb{4}\mywedge \sum \bU{p,q,r}{w-2}  (  - \MC{1}{p-1,q,r+1} (r+1)
+ \MC{2}{p-1,q+1,r} \alpha (q+1) 
 - \MC{3}{p,q,r} ((1 +\alpha) (p+1) + q + \alpha r) )
\end{align*}
\kmcomment{
Assume that \(\pdel(L +  M) = 0\).  
Then we have
\begin{align}
& \label{G3D2:H:0:a}
(1-\alpha)  \LC{ }{p,q+1,r+1} (q+1)(r+1) \\
=& \   
   \MC{1}{p-1,q,r+1} (r+1)
- \MC{2}{p-1,q+1,r} \alpha (q+1) 
 + \MC{3}{p,q,r} ((1 +\alpha) (p+1) + q + \alpha r) 
 \notag
\end{align}
Assume \( 1-\alpha = 0\). Then from  
\eqref{G3D2:H:0:a}, we get  
\begin{align*}
 0 & =    \MC{1}{p-1,q,r+1} (r+1)
    - \MC{2}{p-1,q+1,r} (q+1) 
    + \MC{3}{p,q,r} ( 2+ w + p )
\end{align*}
Thus, \((\MC{3}{p,q,r})\) are completely determined by 
 \(( \MC{1}{p,q,r})\) and \((\MC{2}{p,q,r}) \), and the freedom of 
\((\LC{}{a,b,c})\) and \((\MC{i}{p,q,r})\) is \( \binom{w+2}{2} + 2
\binom{w}{2}\). 
 
When \( 1-\alpha \ne 0\),  \eqref{G3D2:H:0:a} shows that  
\( ( \LC{}{a,b,c} )\) with \(bc >0\)  are determined by 
\( ( \MC{i}{p,q,r} )\) (\(i=1,2,3\). Thus, the
freedom of 
\( ( \LC{}{a,b,c} )\) and 
\( ( \MC{i}{p,q,r} )\) (\(i=1,2,3\)) is \( 2(w+1)-1 + 3 \binom{w}{2}
= \binom{w+2}{2} + 2 \binom{w}{2} \).

\textcolor{red}{Revised} on July 21, 2020 by focusing the rank, the number of generators of linear equations: 
}
The linear equation \( \pdel (L+M)=0\) is generated by 
\begin{align}  \label{G3D2:H:00:a}
& \quad (\alpha -1)  \LC{ }{p,q+1,r+1} (q+1)(r+1) \\
& +    \MC{1}{p-1,q,r+1} (r+1)
- \MC{2}{p-1,q+1,r} \alpha (q+1) 
 + \MC{3}{p,q,r} ((1 +\alpha) (p+1) + q + \alpha r) 
 \notag
\end{align}
If \( w <2\), then generators are 
\( (\alpha -1)  \LC{ }{p,q+1,r+1} (q+1)(r+1) \),  and we get none. Thus, the
rank is 0.

Assume that \( w \geqq 2\). 
If \( \alpha \ne 1\), then we have  \( \binom{w}{2}\) linearly independent
generators 
with leading term \( \LC{}{p,q+1,r+1}\), and the rank is \(\binom{w}{2}\). 
If \( \alpha = 1\), then generators are 
\( \MC{1}{p-1,q,r+1} (r+1) - \MC{2}{p-1,q+1,r} \alpha (q+1) + \MC{3}{p,q,r}
( p + w ) \) with 
leading term \( \MC{3}{p,q,r}\), and the rank is \(\binom{w}{2} \) 
and the kernel dimension is \( 2\binom{w}{2}+ \binom{w+2}{2}\) for two cases. 

\paragraph{The space of cycles in \(\wtedC{w-1}{w}\):} 
Generators of this chain space are \( \zb{4}\mywedge \bU{A}{w-2}\). 

Since 
\(\pdel (\zb{4}\mywedge \bU{A}{w-2}) = 
- \zb{4}\mywedge \pdel(\bU{A}{w-2}) + \SbtES{\zb{4}}{ \bU{A}{w-2}} 
= - \zb{4} \mywedge \zb{4} \wedge  (1-\alpha ) b c \bU{a,b-1,c-1}{w-2} 
+ 0 = 0 \), the boundary operator is trivial and the kernel dimension is 
\( \binom{w}{2}\).

\renewcommand{\arraystretch}{1.3}

\paragraph{Final table of chain complex }%
\begin{thm} \label{thm:G3D2} \ 

\begin{center}
\(
\begin{array}{c|*{5}{c}}
\text{weight}= w > 0 & w-1 & w & w+1 & w+2 & w+3\\\hline
\text{SpaceDim} & 
\binom{w}{2} & 
3 \binom{w}{2} + \binom{w+2}{2} & 
3\binom{w}{2} + 3\binom{w+2}{2} & 
\binom{w}{2} + 3\binom{w+2}{2} & 
\binom{w+2}{2} 
\\
\ker\dim & 
\binom{w}{2} & 
2\binom{w}{2}  +\binom{w+2}{2}  & 
\binom{w}{2}  + 2 \binom{w+2}{2}  & 
\binom{w+2}{2} + \kappa & 
\kappa  \\\hline
\text{Betti} & 0 & 0 & \kappa & 2 \kappa  & \kappa 
\end{array}
\)
\end{center}
where 
\( \kappa = \begin{cases} 1 & \text{if\quad} \alpha = -1 \\
0 & \text{if\quad} \alpha \ne  -1 
\end{cases}
\). 
\end{thm}


\subsection{$\dim[\frakg,\frakg]=3$}
3-dimensional Lie algebras with \(\Sbt{\frakg}{\frakg}=\frakg\) are
completely classified in \cite{MR559927}.  
By an equivalence, called \textit{multiplicativity cogredience}, 
the bracket relations are given by \(
\Sbt{\zb{1}}{\zb{2}}= \zb{3}, 
\Sbt{\zb{1}}{\zb{3}}= -\beta\zb{2}, 
\Sbt{\zb{2}}{\zb{3}}= \alpha \zb{1}\).  
\cite{MR559927} says that if the base field of \(\frakg\) is \(\mR\), then
we may choose \( \alpha = \beta =1\) or \( -\alpha = \beta =1\). 
Those include \(A_{1}\) type by Dynkin diagram. 
Here, we care about two parameters \(\alpha, \beta\) are just non-zero.     

We choose a basis of 2-vectors by \[\ub{1} = \frac{1}{\alpha} \zb{2} \wedge
\zb{3}\;, \; \ub{2} = \frac{-1}{\beta} \zb{1} \wedge \zb{3}\;, \; \ub{3} =
\zb{1} \wedge \zb{2}\] and 3-vector by \( \zb{4} = V= \zb{1} \wedge \zb{2}
\wedge \zb{3} \).   
The multiplication (by the Schouten bracket) tables are the followings: 
\begin{center}
\( \begin{array}[t]{c|*{3}{c}|c|*{3}{c} }
 & \zb{1} & \zb{2} & \zb{3} & \zb{4}& \ub{1} & \ub{2} & \ub{3}  
 \\\hline 
\zb{1} & 0 & \zb{3}  &  -\beta \zb{2} & 0  & 0 & \ub{3} & -\beta\ub{2}   \\
\zb{2} & -\zb{3}  & 0 & \alpha\zb{1}& 0  &  - \ub{3} & 0 & \alpha \ub{1} \\
\zb{3} &  \beta\zb{2} & - \alpha \zb{1} & 0 & 0 & \beta \ub{2}  &  -\alpha \ub{1} & 0  
\end{array}
 \)
\hfil
\( \begin{array}[t]{c|*{3}{c}|c}
   & \ub{1} & \ub{2} & \ub{3} & \zb{4} 
\\\hline 
\ub{1} & \frac{2}{\alpha}\zb{4} & 0 & 0 & 0 \\
\ub{2} & 0 & \frac{2}{\beta}\zb{4} & 0 & 0  \\
\ub{3} & 0 & 0  &  2\zb{4}  & 0\\
\end{array}
 \)
\end{center}

From the tables, we have 
\begin{align*}
\pdel W^{1111} &= 0\;,\quad 
\pdel W^{1110} = 0 \;,\quad 
\\
\pdel W^{1101} &= \zb{3} \mywedge \zb{4} \;,\quad 
\pdel W^{1011} =  -\beta \zb{2} \mywedge \zb{4}\;,\quad 
\pdel W^{0111} =  \alpha \zb{1} \mywedge \zb{4}\;,\quad 
\\
\pdel U^{a,b,c} &= \zb{4} \mywedge (
  \frac{2}{\alpha} \tbinom{a}{2}  U^{a-2, b, c} 
+  \frac{2}{\beta} \tbinom{b}{2}  U^{a, b-2, c} 
+  2 \tbinom{c}{2}  U^{a, b, c-2} 
)
\;. 
\\
\sum_{A} \LC{}{A}\pdel U^{a,b,c} &= \zb{4} \mywedge \sum_{} U^{p,q,r} (
  \frac{2}{\alpha} \tbinom{p+2}{2}  \LC{}{p+2, q, r} 
+  \frac{2}{\beta} \tbinom{q+2}{2}   \LC{}{p, q+2, r} 
+  2 \tbinom{r+2}{2}    \LC{}{p, q, r+2} 
)
\;. 
\\
\SbtES{\zb{1}}{U^{A}} &=  b U^{a,b-1,c+1} - \beta c U^{a,b+1,c-1}\\ 
\SbtES{\zb{2}}{U^{A}} & = - a U^{a-1,b,c+1} + \alpha c U^{a+1,b,c-1} \\
\SbtES{\zb{3}}{U^{A}} &= \beta a U^{a-1,b+1,c}  - \alpha  b U^{a+1,b-1,c}
\\ \sum_{A} \LC{i}{A} \SbtES{\zb{1}}{U^{A}} 
&= \sum_{A}((b+1)\LC{i} {a,b+1,c-1}-\beta (c+1)\LC{i}{a,b-1,c+1}) U^{A}\;,
\\ \sum_{A} \LC{i}{A} \SbtES{\zb{2}}{U^{A}} 
&= \sum_{A}(- (a+1)\LC{i} {a+1,b,c-1}+\alpha (c+1)\LC{i}{a-1,b,c+1}) U^{A}\;,
\\ \sum_{A} \LC{i}{A} \SbtES{\zb{3}}{U^{A}} 
&= \sum_{A}( \beta (a+1)\LC{i} {a+1,b-1,c}-\alpha (b+1)\LC{i}{a-1,b+1,c}) U^{A}\;.
\end{align*}


\paragraph{The space of cycles in \(\wtedC{w+3}{w}\):} 
Take a general chain 
\(L =  \sum_{a+b+c=w} \LC{}{a,b,c} W^{1,1,1,0} \mywedge U^{a,b,c}\) with
unknown scalars \(\LC{}{a,b,c} \). 
\begin{align*} 
\pdel L =&  \pdel ( \sum_{A} \LC{}{A}  W^{1110} \mywedge U^{A})
\\=& \sum \LC{}{a,b,c} ( \pdel   W^{1,1,1,0})  \mywedge U^{A}
-  W^{1110} \sum \LC{}{A} \mywedge \pdel U^{A}
+ \sum_{i=1}^{3} (-1)^{i+1} W^{ \eps{i}=\eps{4}=0 } \sum \LC{ }{A}\SbtES{ \zb{i} } { U^{A
} }
\\ 
=& 
- W^{1111} \mywedge \sum \LC{}{A} 
( 
  \frac{2}{\alpha} \tbinom{a}{2}  U^{a-2, b, c} 
+  \frac{2}{\beta} \tbinom{b}{2}  U^{a, b-2, c} 
+  2 \tbinom{c}{2}  U^{a, b, c-2} 
)
\\& 
+ \zb{2} \mywedge \zb{3} \mywedge \sum \LC{ }{A} 
(  b U^{a,b-1,c+1} - \beta c U^{a,b+1,c-1}
)
\\& 
- \zb{1} \mywedge \zb{3} \mywedge \sum \LC{ }{A} ( 
 - a U^{a-1,b,c+1} + \alpha c U^{a+1,b,c-1} 
)
\\& 
+ \zb{1} \mywedge \zb{2} \mywedge \sum \LC{ }{A}
( 
\beta a U^{a-1,b+1,c}  - \alpha  b U^{a+1,b-1,c}
 )  
\end{align*} 
The linear equation 
\( \pdel L =0\) is generated by the followings:  
\begin{subequations}
\begin{align}
& 
\frac{(p+2)(p+1)}{\alpha }\LC{}{p+2,q,r } 
+\frac{(q+2)(q+1)}{\beta }\LC{}{p,q+2,r } 
+(r+2)(r+1)\LC{}{p,q, r+2} \label{G3:D3:t1:a}
\\
 & 
(b+1) \LC{ }{a,b+1,c-1}  - \beta  (c+1) \LC{ }{a,b-1,c+1} \label{G3:D3:t1:b}
\\& 
 (a+1) \LC{ }{a+1,b,c-1} - \alpha  (c+1) \LC{ }{a-1,b,c+1} \label{G3:D3:t1:c}
\\& 
\beta ( a + 1  ) \LC{}{a+1,b-1,c} -\alpha (b +1) \LC{}{a-1,b+1,c}  
\label{G3:D3:t1:d}
\end{align} 
\end{subequations}


When \(w=0\), then the unknown parameter is only one \(\LC{}{0,0,0}\)
and the above equations give no information. This means  \(\LC{}{0,0,0}\)
is free, and so the kernel dimension is 1.  In fact, \( \pdel ( \zb{1}
\mywedge \zb{2} \mywedge \zb{3}) =0 \). 

When \(w=1\), the unknown parameters  are \(\LC{}{a,b,c}\)
with \(a+b+c=1\). \eqref{G3:D3:t1:a} gives no information, but  
\eqref{G3:D3:t1:d} tells 
\( \LC{}{1,0,0} =0\) , \( \LC{}{0,1,0} =0\).   
\( \LC{}{0,0,1} =0\) comes from 
\eqref{G3:D3:t1:b} by \((a,b,c)=(1,0,0)\). 
Thus, the kernel dimension is 0.  

\kmcomment{
\textcolor{red}{Revised} on July 29, 2020 by focusing the rank, the number
of generators of linear equations:} 
Here we assume \( w\geqq 2\). From  
\eqref{G3:D3:t1:a} $\sim$ \eqref{G3:D3:t1:c}, we have generators of the linear
equations. Here we denote them by 
\begin{align}
G_{0}(p,q,r) & =  
\frac{2}{\alpha } \tbinom{p+2}{2}\LC{}{p+2,q,r } 
+\frac{2}{\beta } \tbinom{q+2}{2}\LC{}{p,q+2,r } 
+ (r+2)(r+1)  \LC{}{p,q, r+2}  \;,
\label{Bb:0:0}
\\   
G_{1}(a,b,c)  &= 
(b+1) \LC{}{a,b+1,c-1}  - \beta (c+1) \LC{}{a,b-1,c+1} \;,
\label{Bb:0:1}
\\ 
G_{2}(a,b,c) &=  
 (a+1)  \LC{}{a+1,b,c-1} - \alpha (c+1)  \LC{}{a-1,b,c+1} \;,
 \label{Bb:0:2}
\\ 
G_{3}(a,b,c) &= 
 \beta (a+1)  \LC{}{a+1,b-1,c} - \alpha (b+1)  \LC{}{a-1,b+1,c} \;.
 \label{Bb:0:3}
\end{align}


\begin{center}
\begin{tikzpicture}[scale=0.7, baseline = (O.base)]
\node (O) at (0,0) {  }; 
\draw[lightgray, ] ( -1,-1  ) grid (1,1); 
\draw node at (0,0) {$\Box$} 
node[] at (0,1) {$\bullet$} 
node[above] at (0,1) {$  (b+1)\LC{}{a,b+1,c-1}$}
node[] at (0,-1) {$\bullet$} 
node[below ] at (0,-1) {$ - \beta (c+1)\LC{}{a,b-1,c+1}$} ;
\end{tikzpicture}
\hfil
\begin{tikzpicture}[scale=0.7, baseline = (O.base)]
\node (O) at (0,0) {  }; 
\draw[lightgray, ] ( -1,-1  ) grid (1,1); 
\draw node at (0,0) {$\Box$} 
node[] at (-1,0) {$\bullet$} 
node[above left] at (0,0) {$ - \alpha (c+1)\LC{}{a-1,b,c+1}$}
node[] at (1,0) {$\bullet$} 
node[below right ] at (0,0) {$ (a+1)\LC{}{a+1,b,c-1}$} ;
\end{tikzpicture}
\hfil
\begin{tikzpicture}[scale=0.7, baseline = (O.base)]
\node (O) at (0,0) {  }; 
\draw[lightgray, ] ( -1,-1  ) grid (1,1); 
\draw node at (0,0) {$\Box$} 
node[] at (-1,1) {$\bullet$} 
node[above ] at (-1,1) {$- \alpha (b+1)\LC{}{a-1,b+1,c}$}
node[] at (1,-1) {$\bullet$} 
node[below ] at (1,-1) {$ \beta (a+1)\LC{}{a+1,b-1,c}$} ;
\end{tikzpicture}
\end{center}

The above pictures suggest properties of \(G_{1}, G_{2}, G_{3}\) are
vertical, horizontal, or slant relationship, where \(\Box\) means the
barycenter.  We call them \DS{} and call \(G_{0} \) \TS{}.  When \(\Box\)
move around the area \(T = \{ (a,b,c) \mid a+b+c=w  \}\), some \(G_{i}(a,b,c)
\) becomes \HS{} near the border of \(T\). For instance,  \(G_{1}(a,0,c) =
\LC{}{a,1,c-1}\). 
By this \HS{} and vertical relations, we see that \( \LC{}{a, 2b+1, c} \)
are \HS{}. By the same idea, we see \( \LC{}{2a+1, b, c} \) and \( \LC{}{a,b,
2c+1 } \) are \HS{}. Take the particular \(\LC{}{0,w-1,1}\) which is \HS{}.
If \(w\) is odd, then \(w-1\) which is the 2nd entry of the \HS{}, and  by
applying vertical relation, \( \LC{}{0, 2b, c}\) are \HS{}. Now applying
horizontal relation, we conclude \(\LC{}{a,b,c} \) are all \HS{} when \(w\)
is odd. So the rank is \( \binom{w+2}{2}\). 

Now assume \(w= 2k\). Since \(\LC{}{a,b,c}\) is \HS{} if one of \(a,b,c\) is
odd, we care about  \(\LC{}{a,b,c}\) of even entries \(a,b,c\).  
Take \begin{align} H_{0} =&  G_{0}( 2p,2q , 2r ) = 
(2p+2)(2p+1) \LC{}{2p+2,2q,2r}  + 
(2q+2)(2q+1) \beta \LC{}{2p,2q+2,2r} \label{H0} \\&   + 
(2r+2)(2r+1) \alpha \LC{}{2p,2q,2r+2}  \notag \\ 
H_{1} =&  G_{1}(2p,1+2q , 1+2r ) = (2+2q) \LC{}{2p,2q+2,2r}- \beta (2+2r)
\LC{}{2p,2q,2r+2} \label{H1}\\
H_{2} =&  G_{2}(1+ 2p, 2q, 1+ 2r ) = (2+2p) \LC{}{2p+2,2q,2r}- \alpha (2+2r)
\LC{}{2p,2q,2r+2} \label{H2} \\
H_{3} =&  G_{3}(1+ 2p, 1+ 2q, 2r ) = \beta (2+2p) \LC{}{2p+2,2q,2r}- \alpha
(2+2q) \LC{}{2p,2q,2r+2} \label{H3}
\end{align}
We see that \( \alpha H_{1} - \beta  H_{2} + H_{3} =0\) and  
\begin{align}
H_{0} =& \frac{2p+1}{\alpha} \left( H_{2} + \alpha (2+2r)\LC{}{2p,2q,2r+2} \right)
 + \frac{2q+1}{\beta} \left( H_{1} + \beta (2+2r)\LC{}{2p,2q,2r+2} \right)
 \\
& +  (2r+2) (2r+1) \alpha \LC{}{2p,2q,2r+2} 
\\
=& \frac{2p+1}{\alpha} H_{2} + \frac{2q+1}{\beta} H_{1} + \alpha (2w-1) (2r+2) \LC{}{2p,2q,2r+2}
\end{align}
This shows \( \LC{}{2p,2q,2r+2} \) is a \HS{}.  We also see that 
\( \LC{}{2p+2,2q,2r} \) is a \HS{} by \eqref{H2} and  
\( \LC{}{2p,2q+2,2r} \) is a \HS{} by \eqref{H3}. Thus, every \(\LC{}{a,b,c}
\) is \HS{} and the rank is \( \binom{w+2}{2}\).  

\bigskip

\paragraph{The space of cycles in \(\wtedC{w+2}{w}\):} 
We study the space of cycles in \(\wtedC{w+2}{w}\). 
\kmcomment{where  
the generators are  \( W[\eps{i}=\eps{4}=0]
\mywedge \bU{a,b,c}{w}\)
or \( W^{1,1,1,1} \mywedge \bU{p,q,r}{w-2}\).   
}
If \(w=0\), \( \zb{i} \mywedge \zb{j}\) (\(i<j\leqq 3\)) are a basis and 
since \(\pdel( \zb{1} \mywedge \zb{2} )= \zb{3}\),  
\(\pdel( \zb{1} \mywedge \zb{3} )= - \beta \zb{2}\),  and 
\(\pdel( \zb{2} \mywedge \zb{3} )=  \alpha \zb{1}\),   the rank is 3 and the
kernel dimension is 0. 
If \(w=1\), \( W[\eps{i}=\eps{4}=0] \mywedge \ub{j} \) are a basis of the
chain space.  Take a general chain 
\( X = \sum_{i,j} \LC{i}{j}  W[\eps{i}=\eps{4}=0] \mywedge \ub{j} \).   
\begin{align*}
\pdel X =& \sum_{j} \LC{1}{j} ( \alpha \zb{1} \mywedge \ub{j} 
+ \zb{3} \mywedge \Sbt{\zb{2}} {\ub{j}}
- \zb{2} \mywedge \Sbt{\zb{3}} {\ub{j}}) \\
&+ \sum_{j} \LC{2}{j} ( - \beta \zb{2} \mywedge \ub{j} 
+ \zb{3} \mywedge \Sbt{\zb{1}} {\ub{j}}
- \zb{1} \mywedge \Sbt{\zb{3}} {\ub{j}}) \\
&+ \sum_{j} \LC{3}{j} (  \zb{3} \mywedge \ub{j} 
+ \zb{2} \mywedge \Sbt{\zb{1}} {\ub{j}}
- \zb{1} \mywedge \Sbt{\zb{2}} {\ub{j}}) \\
=& 
+  \LC{1}{1} ( \alpha \zb{1} \mywedge \ub{1} 
+ \zb{3} \mywedge (- {\ub{3}} ) 
- \zb{2} \mywedge (\beta {\ub{2}}) \\
& + \LC{1}{2} ( \alpha \zb{1} \mywedge \ub{2} 
+ \zb{3} \mywedge 0
- \zb{2} \mywedge ( - \alpha  {\ub{1}}) 
+  \LC{1}{3} ( \alpha \zb{1} \mywedge \ub{3} 
+ \zb{3} \mywedge \alpha  {\ub{1}}
- \zb{2} \mywedge 0 ) \\
&+  \LC{2}{1} ( - \beta \zb{2} \mywedge \ub{1} 
+ \zb{3} \mywedge  0 
- \zb{1} \mywedge \beta {\ub{2}}) \\
&+  \LC{2}{2} ( - \beta \zb{2} \mywedge \ub{2} 
+ \zb{3} \mywedge  {\ub{3}}
- \zb{1} \mywedge ( - \alpha  {\ub{1}}) ) 
+ \LC{2}{3} ( - \beta \zb{2} \mywedge \ub{3} 
+ \zb{3} \mywedge ( - \beta  {\ub{2}}) 
- \zb{1} \mywedge 0 ) \\
&+ \LC{3}{1} (  \zb{3} \mywedge \ub{1} 
+ \zb{2} \mywedge 0 
- \zb{1} \mywedge ( -  {\ub{3}})) \\
&+ \LC{3}{2} (  \zb{3} \mywedge \ub{2} 
+ \zb{2} \mywedge  {\ub{3}}
- \zb{1} \mywedge 0 ) 
+ \LC{3}{3} (  \zb{3} \mywedge \ub{3} 
+ \zb{2} \mywedge ( - \beta  {\ub{2}} ) 
- \zb{1} \mywedge \alpha   {\ub{1}}) \\
\end{align*}
\(\pdel X=0\) is defined by 6 linearly independent polynomials
\( \LC{1}{1} +  \LC{2}{2} -  \LC{3}{3}\),   
\( \LC{1}{1} +  \LC{2}{2} +  \LC{3}{3}\),   
\( - \LC{1}{1} +  \LC{2}{2} +  \LC{3}{3}\),   
\( \alpha  \LC{1}{2}  - \beta  \LC{2}{1} \),   
\( \alpha  \LC{1}{3}  +   \LC{3}{1} \),   
\( - \beta  \LC{2}{3}  +   \LC{3}{2} \). Thus the rank is 6 and the kernel
dimension -s 3.    

When \( w>1\), 
take a general chain \(L+M\) with 
\(L= \sum \LC{i}{a,b,c}  W[{\eps{i}=\eps{4}=0}] \mywedge \bU{a,b,c}{w}
\) and \( M= 
 \sum \MC{} {p,q,r}  W^{1111} \mywedge \bU{p,q,r}{w-2}
\) with unknown scalars \(\LC{i}{a,b,c}, \MC{}{p,q,r}\).   

Then 
\begin{align*}
& \pdel (L+ M) \\ 
= & 
 \sum \LC{i}{A}( \pdel( W[{\eps{i}=\eps{4}=0}]) \mywedge U^{A}
+  W[{\eps{i}=\eps{4}=0}] \mywedge  \pdel U^{A} )
+ \SbtES{
W[{\eps{i}=\eps{4}=0}] } { U^{A} } 
\\ & 
+ 0 + 0 + \sum \MC{} {P}  \SbtES{ W^{1111}}{ U^{P}}  
\\=&    
 \alpha  \zb{1} \mywedge\sum \LC{1}{A} U^{A}
- \beta  \zb{2} \mywedge\sum \LC{2}{A} U^{A}
+  \zb{3} \mywedge\sum \LC{3}{A} U^{A}
\\& 
+ W^{0111}\mywedge\sum\LC{1}{A}(
  \frac{2}{\alpha} \tbinom{a}{2}  U^{a-2, b, c} 
+  \frac{2}{\beta} \tbinom{b}{2}  U^{a, b-2, c} 
+  2 \tbinom{c}{2}  U^{a, b, c-2} 
) 
\\& 
+ W^{1011}\mywedge\sum\LC{2}{A}(
\frac{2}{\alpha} \tbinom{a}{2}  U^{a-2, b, c} 
+  \frac{2}{\beta} \tbinom{b}{2}  U^{a, b-2, c} 
+  2 \tbinom{c}{2}  U^{a, b, c-2} 
) 
\\& 
+ W^{1101}\mywedge\sum\LC{3}{A}(
  \frac{2}{\alpha} \tbinom{a}{2}  U^{a-2, b, c} 
+  \frac{2}{\beta} \tbinom{b}{2}  U^{a, b-2, c} 
+  2 \tbinom{c}{2}  U^{a, b, c-2} 
) 
\\& 
+ \sum \LC{1}{A} ( - \zb{2} \SbtES{ \zb{3} }{ U^{A} }
                   + \zb{3} \SbtES{ \zb{2} }{ U^{A} } )
+ \sum \LC{2}{A} ( - \zb{1} \SbtES{ \zb{3} }{ U^{A} }
                   + \zb{3} \SbtES{ \zb{1} }{ U^{A} } )
                   \\& 
+ \sum \LC{3}{A} ( - \zb{1} \SbtES{ \zb{2} }{ U^{A} }
                   + \zb{2} \SbtES{ \zb{1} }{ U^{A} } )
\\&
+ W^{0111} \mywedge \MC{}{P} \SbtES{ \zb{1} }{ U^{P} } 
- W^{1011} \mywedge \MC{}{P} \SbtES{ \zb{2} }{ U^{P} } 
+ W^{1101} \mywedge \MC{}{P} \SbtES{ \zb{3} }{ U^{P} } 
\\
= & 
+  \zb{1} \mywedge\sum ( \alpha \LC{1}{A} U^{A}
 -  \LC{2}{A} \SbtES{ \zb{3} }{ U^{A} }
 -  \LC{3}{A} \SbtES{ \zb{2} }{ U^{A} }
 )
\\& 
+  \zb{2} \mywedge\sum (- \beta \LC{2}{A} U^{A}
  -  \LC{1}{A}\SbtES{ \zb{3} }{ U^{A} }
 +  \LC{3}{A} \SbtES{ \zb{1} }{ U^{A} } 
 )
\\& 
+  \zb{3} \mywedge\sum ( \LC{3}{A} U^{A}
  +  \LC{1}{A}\SbtES{ \zb{2} }{ U^{A} } 
 +  \LC{2}{A} \SbtES{ \zb{1} }{ U^{A} } 
 )
\\& 
+ W^{0111} \mywedge (\sum  \MC{}{P} \SbtES{\zb{1}}{U^{P}} 
+\sum \myL{1}{p,q,r} U^{p,q,r})  
\\& 
+ W^{1011} \mywedge \sum (- \MC{}{P} \SbtES{\zb{2}}{U^{P}} 
+\sum \myL{2}{p,q,r} U^{p,q,r}) 
\\& 
+ W^{1101} \mywedge\sum ( \MC{}{P} \SbtES{\zb{3}}{U^{P}} 
+\sum  \myL{3}{p,q,r} U^{p,q,r})  
\\\shortintertext{where}
\myL{i}{p,q,r} &= 
 \frac{(p+2)(p+1)}{\alpha} \LC{i}{p+2,q,r}
 +\frac{(q+2)(q+1)}{\beta} \LC{i}{p,q+2,r}
 +(r+2)(r+1) \LC{i}{p,q,r+2} \; . 
\end{align*}
Then \( \pdel (L+M) =0\) is defined by the following polynomials of unknown
\(\LC{i}{A}\) and \(\MC{}{P}\): 
\begin{align}
\myE{1}{a,b,c} =& 
\alpha \LC{1}{a,b,c}  
 - \beta (a+1) \LC{2}{a+1,b-1,c} + \alpha (b+1) \LC{2}{a-1,b+1,c} 
 \label{G3:D3:type1:CwP2:a}
 \\&\quad 
 + (a+1) \LC{3} {a+1,b,c-1} - \alpha (c+1) \LC{3}{a-1,b,c+1} 
 \notag
\\
\myE{2}{a,b,c} = & 
- \beta  \LC{2}{a,b,c} 
 -  \beta (a+1) \LC{1}{a+1,b-1,c} + \alpha (b+1) \LC{1}{a-1,b+1,c} 
\label{G3:D3:type1:CwP2:b}
 \\&\quad 
 +  (b+1) \LC{3}{a,b+1,c-1} - \beta (c+1) \LC{3}{a,b-1,c+1}
 \notag
\\
\myE{3}{a,b,c} = &
 \LC{3}{a,b,c} 
 - (a+1) \LC{1} {a+1,b,c-1} + \alpha (c+1) \LC{1}{a-1,b,c+1} 
 \label{G3:D3:type1:CwP2:c} 
 \\&\quad 
 +  (b+1) \LC{2}{a,b+1,c-1} - \beta (c+1) \LC{2}{a,b-1,c+1}
 \notag
\\
\myE{4}{p,q,r} = & 
( q+1) \MC{}{p,q+1,r-1}  - \beta (r+1) \MC{}{p,q-1,r+1}  
+ \myL{1}{p,q,r} 
\label{G3:D3:type1:CwP2:d} 
\\
\myE{5}{p,q,r} = & 
(p+1) \MC{}{p+1,q,r-1}  - \alpha (r+1) \MC{}{p-1,q,r+1}
+ \myL{2}{p,q,r}
 \label{G3:D3:type1:CwP2:e}  
\\
\myE{6}{p,q,r} = &
  \beta (p+1)\MC{}{p+1,q-1,r} 
  - \alpha (q+1)\MC{}{p-1,q+1,r} 
+ \myL{3}{p,q,r}
 \label{G3:D3:type1:CwP2:f} 
\end{align}


\subsubsection{Rank of \( \MC{}{p,q,r}\) }
We see that \( \myE{6}{p,q,r} - \beta \myE{5}{p,q-1,r+1} + \alpha \myE{4}{p-1,
q,r+1} \) does not have \( \MC{}{P}\) as follows:
\begin{align}
& 
+p (p+1) \LC{1}{ p+1,q,r+1 }
+\alpha (r+3) (r+2) \LC{1}{ p-1,q,r+3 }
+\alpha (q +2) (q+1)/\beta \LC{1}{ p-1,q+2,r+1 }
\label{side:effect:mu}
\\& 
-q (q+1) \LC{2}{ p,q+1,r+1 }
-\beta (r+3) (r+2) \LC{2}{ p,q-1,r+3 }
-\beta (p+2) (p+1)/\alpha \LC{2}{ p+2,q-1,r+1 }
\notag
\\&
+(r+2) (r+1) \LC{3} { p,q,r+2 }
+(p+2) (p+1)/\alpha \LC{3}{ p+2,q,r }
+(q+2) (q+1)/\beta \LC{3}{ p,q+2,r } \; .
\notag
\end{align}
In order to know the rank concerning to \(\MC{}{p,q,r}\), we  handle
one of three, \( \myF{6}{} =  \myE{6}{}\negmedspace \mod\LC{i}{A}\),
i.e.,   
\(\myF{6}{p,q,r}=\beta (p+1)\MC{}{p+1,q-1,r}-\alpha (q+1)\MC{}{p-1, q+1,r}\)
. \kmcomment{
and \(\myF{5}{} =  \myE{5}{}\negthickspace \mod\LC{i}{A} \), i.e., \(\myF{5}
{p,q,r} = (p+1) \MC{}{p+1,q,r-1}  - \alpha (r+1) \MC{}{p-1,q,r+1} \).  
}
Among \( \myF{6}{p,q,r}\) with \(p+q+r = w-2\), first fix \(r\) and put
\(H=w-2 - r\) and study of linearly independence of \(\Hori{w-2}{r} = \{
\myF{6}{p,q,r} \mid p+q = H=w-2-r\}\) as polynomials of  \(\MC{}{P}\).
Basically, \( \myF{6}{p, q,r}\) looks like  slant double stars with distance
2,  we define 


\begin{align*} 
OD^{[r]}_{2i+1} & = \myF{6}{1,H-1,r} \wedge 
\myF{6}{3,H-3,r} \wedge \cdots \wedge \myF{6}{2i+1,H-1-2i,r} 
\quad\mathtt{, where} \quad   2i+1 \leqq H\; ,   
\\
&=  \mathop{\wedge}_{i=0}^{2i+1\leqq H} ( -\alpha( H-2i ) \MC{}{2i,H-2i,r} + \beta
(2+2i) \MC{}{2i+2,H-2i-2,r} )
\\\shortintertext{and}
EV^{[r]}_{2i} & = \myF{6}{0,H,r} \wedge 
\myF{6}{2,H-2,r} \wedge \cdots \wedge \myF{6}{2i,H-2i,r} \quad\mathtt{where}
\quad  2i \leqq H\;. \\\shortintertext{For instance,} 
EV^{[r]}_{0} & = \myF{6}{0,H,r} = 
\beta \MC{}{1,H-1,r}-\alpha (H+1)\MC{}{-1,H+1,r} = \beta \MC{}{1,H-1,r}\;, 
\\
EV^{[r]}_{2} & = EV^{[r]}_{0} \wedge \myF{6}{2,H-2,r} =  \beta \MC{}{1,H-1,r}\wedge  
(\beta 3 \MC{}{3,H-3,r}-\alpha (q+1)\MC{}{1,H-1,r}) 
\\& 
= 3 \beta^{2} \MC{}{1,H-1,r}\wedge  \MC{}{3,H-3,r}
\\ & \vdots \\
 EV^{[r]}_{2k} &= (2k+1)!! \beta^{k+1} 
\MC{}{1,H-1,r}\wedge  \MC{}{3,H-3,r} \wedge \cdots \wedge \MC{}{2k+1,
H-2k-1,r}
\\\shortintertext{If \(H= 2  h\) even, then }
 EV^{[r]}_{2h} &= EV^{[r]}_{2h-2}  \wedge \myF{6}{2h,0,r} \\
& = (2h-1)!! \beta^{h} 
\MC{}{1,H-1,r}\wedge \wedge \cdots \wedge \MC{}{2h-1,1,r}
\wedge
( 0 -\alpha \MC{}{2h-1,1,r}) 
\\& = 0 \; . 
\end{align*}

Thus, we have the next lemma. 
\begin{lemma}
\label{lemma:slant:mu}
If \(H=w-2-r= \) odd, then 
the set \(\{ \myF{6}{p,q,r} \mid p+q= H=w-2 - r\}\) consists of linearly
independent polynomials. 
\kmcomment{
and so  the set \(\{ \myE{6}{p,q,r} \mid p+q= H=w-2 - r\}\) consists of
linearly independent polynomials.}

If \(H=w-2-r= \) even, then 
the set \(\{ \myF{6}{p,q,r} \mid p+q= H=w-2 - r\}\) is not linearly
independent, but   
the set \(\{ \myF{6}{p,q,r} \mid p+q= H=w-2 - r\} \setminus 
\{ \myF{6}{H,0,r} \} \) consists of linearly independent polynomials. 
\kmcomment{but    
the set \(\{ \myE{6}{p,q,r} \mid p+q= H=w-2 - r\} \)
 is linearly independent.}    
\end{lemma}
\kmcomment{
\textbf{Proof:} 
What is left to be proven is the last statement of Lemma. 
We recall \[ \myE{6}{p,q,r} = \myF{6}{p,q,r} 
+ \frac{(p+2)(p+1)} {\alpha} \LC{3}{p+2,q,r} 
+ \frac{(q+2)(q+1)} {\beta} \LC{3}{p,q+2,r} 
+ (r+2)(r+1) \LC{3}{p,q,r+2} \]
and simply prove that when \(r\) is fixed, they are linearly independent
focusing \(\LC{3}{p,q,r}\).  
\kmqed 
}
As a direct corollary, we see that 
\begin{kmCor} 
The dimension of the space spanned by \( \myE{4}{}, \myE{5}{}, \myE{6}{}\)  with leading term
\( \MC{}{P}\) is \(\binom{w}{2}\) if \(w\) is odd, and is \(\binom{w}{2}-1\)
when \(w\) is even.       
\end{kmCor} 

\subsubsection{Rank of \(\LC{i}{a,b,c}\) }

We define \(\myWE{1}{a,b,c}\) symbolically by  \(\myE{1}{a,b,c}/\myE{3}{A}\)
and \(\myWE{2}{a,b, c}\) by \( \myE{2}{a,b,c}/\myE{3}{A}\)  because  all
\( \LC{3}{a',b',c'}\) are completely determined by \( \LC{1}{A}, \LC{2}{A}\)
in \eqref{G3:D3:type1:CwP2:c}. 
\(\myWE{1}{}, \myWE{2}{}\) are linear polynomials of only  
\( \LC{1}{A}, \LC{2}{A}\) as follows: 
\begin{align}
\myWE{1}{a,b,c} &= 
+\alpha ^2 (c+2) (c+1)  \LC{1}{a-2,b,c+2}
- \alpha  (2 a c+a+c-1)  \LC{1}{a,b,c}
+ (a+2) (a+1)  \LC{1}{a+2,b,c-2}  \label{CwP2:EE1}
\\& 
+\alpha  (c+2) (b+1)  \LC{2}{a-1,b+1,c}
-(b+1) (a+1)  \LC{2}{a+1,b+1,c-2} \notag
\\& 
-\alpha  \beta  (c+2) (c +1)  \LC{2}{a-1,b-1,c+2}
+\beta  (c- 1) (a+1)  \LC{2}{a+1,b-1,c} \notag
\\
\myWE{2}{a,b,c} &= 
-\alpha  (c-1) (b+1)  \LC{1}{a-1,b+1,c}
+ (a+1)(b+1)   \LC{1}{a+1,b+1,c-2} \label{CwP2:EE2}
\\& 
+\alpha  \beta  (c+2) (c+1)  \LC{1}{a-1,b-1,c+2 }
-\beta  (a+1) (c+2) \LC{1}{a+1,b-1,c} \notag
\\& 
- (b+2) (b+1)  \LC{2}{a,b+2,c-2}
+\beta  (2 b c+ b + c-1)  \LC{2}{a,b,c}
-\beta ^2 (c+2) (c+1)  \LC{2}{a,b-2,c+2} \notag
\end{align}
In \eqref{CwP2:EE1} or 
\eqref{CwP2:EE2}, almost all coefficients of each term are non-zero, but we
see \(c-1\) , \( 2ac + a+c-1\) or \(2bc+b+c-1\).  
We get special cases as below:
\begin{align}
\frac{1}{\alpha} \myWE{1}{a,b,1} &= 
\frac{1}{\beta} \myWE{2}{a-1,b+1,1}\;,  \\ 
 \myWE{1}{0,w-1,1} &= 0\;, \quad \myWE{2}{w-1,0,1} = 0 \; . 
\end{align}

\kmcomment{
For a fix weight w, we have linear polynomials \( \myWE{1}{a,b,c}\) of \( \{
\LC{1}{a', b',c'} \mid a'+b'+c' =w \}\)  for each \((a,b,c)\) of \(\mN^{3}\)
with \(a+b+c= w\).    
}

We saw  \( \myE{6}{p,q,r} - \beta \myE{5}{p,q-1,r+1} + \alpha \myE{4}{p-1,
q,r+1} \), (say \(G\)) does not have \( \MC{}{P}\) in \eqref{side:effect:mu}. 
We define \(\widetilde{G} = G/\myE{3}{A}\). Then we have  
\begin{align*}
\widetilde{G} = & 
+(p+3) (p+2) (p+1)/\alpha \LC{1}{p+3,q,r-1}
+(q+2) (q+1) (p+1)/\beta \LC{1}{p+ 1,q+2,r-1}
\\& 
-r (p+1) (p-r-1)  \LC{1}{p+1,q,r+1}
-\alpha r (r+3) (r+2) \LC{1}{p-1,q, r+3}
\\& 
-\alpha r (q+2) (q+1)/\beta \LC{1}{p-1,q+2,r+1}
+r (q+1) (q-r-1)  \LC{2}{p,q+ 1,r+1}
\\& 
+\beta r (r+3) (r+2)  \LC{2}{p,q-1,r+3}
+\beta r (p+2) (p+1)/\alpha  \LC{2}{p+2, q-1,r+1}
\\& 
-(p+2) (p+1) (q+1)/\alpha \LC{2}{p+2,q+1,r-1}
-(q+3) (q+2) (q+1)/ \beta \LC{2}{p,q+3,r-1}
\end{align*}

Now we define \(\myF{1}{a,b,c}\) symbolically by  \(\myWE{1}{a,b,c}/\LC{2}
{A}\).  Since \(\myF{1}{a,b,c}\) looks like horizontal 3 stars with 2 step distance
in its shape, we study of linear independence of the subset 
\( \{ \myF{1} {a,b,c} \mid a+c = w - b \}\) and 
\( \{ \myWE{1} {a,b,c} \mid a+c = w - b \}\) 
for each \(b\)  ( \(0 \leqq b \leqq  w\) ) by the same discussion 
in the previous subsection, there are  
 slant double stars with 2 step distance in its
shape.

When  \(a+b+c=w\), we fix  \(b\) and abbreviate  \(w-b\) as
\(H\),  \(c\) is determined by 
\(c = H-a \), we use the following notations 
\[ \myF{1}{a,b,c} = 
\myF{1,w}{a,b,c} = \myF{1,w,b}{a}= \myF{1,w}{a}= \myF{1}{a}
\;, \;  
\LC{1}{a,b,c} = 
\Lc{w}{a,b,c} = \Lc{w,b}{a} = \Lc{w}{a} = \Lc{}{a} 
\] unless there is confusion.  
The following are  matrices of \(\myF{1}{}\) with respect \( \LC{1}{}\)  in
which the coefficients of simultaneous linear equations are arranged. Since
the size is large, we write them separately with an even number of
subscripts. They are tri-diagonal matrices.  

Case of $\mikH=w-b = 1+2h=$odd :
\setlength{\arraycolsep}{1.5pt}
\begin{center}
\(
\begin{array} {c | *{14}{c} } 
\hline
 & \Lc{w,b}{0} & \Lc{w,b}{2} & \Lc{w,b}{4} &&& & & \Lc{w,b}{2h-4} & \Lc{w,b}{2h-2} & \Lc{w,b}{2h} \\
 \hline
\myF{1,w,b}{0} &  -\alpha (H-1) & 2 & 0 \\
\myF{1,w,b}{2} & \alpha^2 H(H-1)& * &  3 \cdot 4  \\
\\
\myF{1,w,b}{a} & &&&  \alpha^2 (1+c)(2+c) & * &  (1+a) (2+a)  \\
\\
\myF{1,w,b}{2h-2} &&& &&&&& 20 \alpha^2  & * &  (H-1) (H-2)  \\
\myF{1,w,b}{2h} &&& &&& &&  0  & 6\alpha^{2} &  - 3\alpha (H-1)  \\
\hline
\end{array}
\)
\end{center}

\bigskip

\begin{center}
\(
\begin{array} {c | *{14}{c} } 
\hline
 & \Lc{w,b}{1} & \Lc{w,b}{3} & \Lc{w,b}{5} &&& & & \Lc{w,b}{2h-3} & \Lc{w,b}{2h-1} & \Lc{w,b}{2h+1} \\
 \hline
\myF{1,w,b}{1} &  -3\alpha (H-1) & 2\cdot 3 & 0 \\
\myF{1,w,b}{3} & \alpha^2 (H-1)(H-2)& * &  4 \cdot 5  \\
\\
\myF{1,w,b}{a} & &&&  \alpha^2 (1+c)(2+c) & * &  (1+a) (2+a)  \\
\\
\myF{1,w,b}{2h-1} &&& &&&&&  12 \alpha^2  & * &  H (H-1)  \\
\myF{1,w,b}{2h+1} &&& &&& &&  0  & 2\alpha^{2} &  - \alpha (H-1)  \\
\hline
\end{array}
\)
\end{center}
where 
\(* = \alpha \left( 1- (a+c+2 a c) \right) \).  

From the shape of matrices, their rank is \(h\) or \(h+1\), and we want to
show the last line is not linearly independent. 

\bigskip

Method of deforming \TS{} to \DS{}: This discussion works for any parity
\(\mikH\).  
Assume that we have \TS{}
\(\myTF{1, w,b}{a} = \mykappa{a} \Lc{w,b}{a} + (a+1)
(a+2) \Lc{w,b} {a+2} \) with \( \mykappa{a} \ne 0 \)  
(from \DS{}
\( \myF{1,w,b} {a }\)).     
We already have 
\(\myTF{1, w,b}{0} = \myF{1, w,b}{0} \) and 
\(\myTF{1, w,b}{1} = \myF{1, w,b}{1} \) with  
\( \mykappa{0} = - \alpha (H-1) \),
\( \mykappa{1} = -3 \alpha (H-1) \).  

Define a \DS{} by 
\begin{align} 
 \myTF{1, w,b}{a+2} & = \myF{1,w,b}{a+2} - \frac{ \alpha^2 (H-a)(H-1-a) }
 {\mykappa {a}} \myTF{1, w,b}{a} 
 \\
& = \mykappa{a+2} \Lc{w,b}{a+2} + (a+3)(a+4) \Lc{w,b}{a
 +4}
\\ 
 \mykappa{a+2} & = - \alpha( 2H a + 5 H - 2 a^2 - 8 a -9 )
 - \frac{ \alpha^2 (H-a)(H-1-a) }{\mykappa{a}} (a+1) (a+2)
\\\shortintertext{
Inductively, 
starting from 0 by step 2, we have } \mykappa{2a} &= - \alpha  (2a+1) (H-1-2a) 
\\\shortintertext{starting from 1 by step2 similarly, we have }  
\mykappa{2a+1} & = -  \alpha (2 + (2a+1) ) (H-(2a+1)) 
\end{align}

\begin{center}
\setlength{\arraycolsep}{5pt}
\(
\begin{array} {c | *{14}{c} } 
\hline
 & \Lc{w,b}{0} & \Lc{w,b}{2} & \Lc{w,b}{4} &&& & & \Lc{w,b}{2h-4} & \Lc{w,b}{2h-2} & \Lc{w,b}{2h} \\
 \hline
\myTF{1,w,b}{0} &  -\alpha (H-1) & 2 & 0 \\
\myTF{1,w,b}{2} &  0 & \mykappa{2} &  3 \cdot 4  \\
\\
\myTF{1,w,b}{a} & &&&  0  & \mykappa{a} &  (1+a) (2+a)  \\
\\
\myTF{1,w,b}{2h-2} &&& &&&&& 0   & \mykappa{2h-2} &  (H-1) (H-2)  \\
\myTF{1,w,b}{2h} &&& &&& &&  0  & 6\alpha^{2} &  - 3\alpha (H-1)  \\
\hline
\end{array}
\)


\(
\begin{array} {c | *{14}{c} } 
\hline
 & \Lc{w,b}{1} & \Lc{w,b}{3} & \Lc{w,b}{5} &&& & & \Lc{w,b}{2h-3} & \Lc{w,b}{2h-1} & \Lc{w,b}{2h+1} \\
 \hline
\myTF{1,w,b}{1} &  -3\alpha (H-1) & 2\cdot 3 & 0 \\
\myTF{1,w,b}{3} & 0& \mykappa{3} &  4 \cdot 5  \\
\\
\myTF{1,w,b}{a} & &&&  0  &  \mykappa{a} &  (1+a) (2+a)  \\
\\
\myTF{1,w,b}{2h-1} &&& &&&&&  0  & \mykappa{2h-1} &  H (H-1)  \\
\myTF{1,w,b}{2h+1} &&& &&& &&  0  & 2\alpha^{2} &  - \alpha (H-1)  \\
\hline
\end{array}
\)
\end{center}

\bigskip

\(\mikH=w-b= 2h+1\)  
\[ \mykappa{2h-2} = -2 (2h-1) \alpha\;,\quad 
 \mykappa{2h-1} = -2 (2h+1) \alpha\;,\quad 
 \begin{vmatrix} \mykappa{2h-2} & (H-1) (H-2) \\
6\alpha^2 & -3\alpha (H-1) \end{vmatrix}
= \begin{vmatrix} -2(2h-1) \alpha & 2h (2h-1) \\
6\alpha^2 & -3\alpha 2h  \end{vmatrix} = 0\] \[ 
 \begin{vmatrix} \mykappa{2h-1} & H (H-1) \\
2\alpha^2 & -\alpha (H-1) \end{vmatrix}
= \begin{vmatrix} -2(2h+1) \alpha & 2h (2h+1) \\
2\alpha^2 & -\alpha 2h  \end{vmatrix} = 0\]
Thus, 
\begin{align*}
& \{ \myTF{1}{2a,b,H-2a} \mid 0\leqq 2a \leqq H-3\} \quad \text{are linearly
independent, and }\\
& \{ \myTF{1}{2a,b,H-2a} \mid 0\leqq 2a \leqq H-1\} \quad \text{are not linearly
independent. }\\
& \{ \myTF{1}{2a+1,b,H-1-2a} \mid 1\leqq 2a+1 \leqq H-2\} \quad \text{are linearly
independent, and }\\
& \{ \myTF{1}{2a+1,b,H-1-2a} \mid 1\leqq 2a+1 \leqq H\} \quad \text{are not linearly
independent. }
\end{align*}

\medskip

Case of $\mikH=w-b = 2h=$even : 
\( \#\{ \myF{1,w,b}{2a} \} = \# \{ \Lc{w,b}{2a} \} = h+1 \ne h = 
 \#\{ \myF{1,w,b}{2a+1} \} = \# \{ \Lc{w,b}{2a+1} \}\). This is different
from the case of \(\mikH=\) odd.  

\begin{center}
\setlength{\arraycolsep}{1.5pt}
\(
\begin{array} {c | *{14}{c} } 
\hline
 & \Lc{w,b}{0} & \Lc{w,b}{2} & \Lc{w,b}{4} &&& & & \Lc{w,b}{2h-4} & \Lc{w,b}{2h-2} & \Lc{w,b}{2h} \\
 \hline
\myF{1,w,b}{0} &  -\alpha (H-1) & 2 & 0 \\
\myF{1,w,b}{2} & \alpha^2 H(H-1)& * &  12  \\
\\
\myF{1,w,b}{a} & &&&  \alpha^2 (1+c)(2+c) & * &  (1+a) (2+a)  \\
\\
\myF{1,w,b}{2h-2} &&& &&&&& 12 \alpha^2  & * &  H (H-1)  \\
\myF{1,w,b}{2h} &&& &&& &&  0  & 2\alpha^{2} &  - \alpha (H-1)  \\
\hline
\end{array}
\)

\bigskip

\(
\begin{array} {c | *{9}{c} } 
\hline
 & \Lc{w,b}{1} & \Lc{w,b}{3} & \Lc{w,b}{5} && & & \Lc{w,b}{2h-5} & \Lc{w,b}{2h-3} & \Lc{w,b}{2h-1} \\
 \hline
\myF{1,w,b}{1} &  -3\alpha (H-1) & 6 & 0 \\
\myF{1,w,b}{3} & \alpha^2 (H-1)(H-2)& * &  20 \\
\\
\myF{1,w,b}{a} & &&&  \alpha^2 (1+c)(2+c) & * &  (1+a) (2+a)  \\
\\
\myF{1,w,b}{2h-3} &&& &&&&  20 \alpha^2  & * &  (H-1) (H-2) \\
\myF{1,w,b}{2h-1} &&& &&&&  0  & 6\alpha^{2} &  -3 \alpha (H-1)  \\
\hline
\end{array}
\)
\end{center}

\bigskip 
   
\begin{center}
\setlength{\arraycolsep}{5pt}
\(
\begin{array} {c | *{14}{c} } 
\hline
 & \Lc{w,b}{0} & \Lc{w,b}{2} & \Lc{w,b}{4} &&& & & \Lc{w,b}{2h-4} & \Lc{w,b}{2h-2} & \Lc{w,b}{2h} \\
 \hline
\myTF{1,w,b}{0} &  -\alpha (H-1) & 2 & 0 \\
\myTF{1,w,b}{2} &  0 & \mykappa{2} &  12  \\
\\
\myTF{1,w,b}{a} & &&&  0  &  \mykappa{a} &  (1+a) (2+a)  \\
\\
\myTF{1,w,b}{2h-2} &&& &&&&& 0  & \mykappa{2h-2}  &  H (H-1)  \\
\myTF{1,w,b}{2h} &&& &&& &&  0  & 2\alpha^{2} &  - \alpha (H-1)  \\
\hline
\end{array}
\)

\bigskip

\(
\begin{array} {c | *{9}{c} } 
\hline
 & \Lc{w,b}{1} & \Lc{w,b}{3} & \Lc{w,b}{5} && & & \Lc{w,b}{2h-5} & \Lc{w,b}{2h-3} & \Lc{w,b}{2h-1} \\
 \hline
\myTF{1,w,b}{1} &  -3\alpha (H-1) & 6 & 0 \\
\myTF{1,w,b}{3} &  0 & \mykappa{3} &  20 \\
\\
\myTF{1,w,b}{a} & &&&   0  & \mykappa{a} &  (1+a) (2+a)  \\
\\
\myTF{1,w,b}{2h-3} &&& &&&&  0   & \mykappa{2h-3} &  (H-1) (H-2) \\
\myTF{1,w,b}{2h-1} &&& &&&&  0  & 6\alpha^{2} &  -3 \alpha (H-1)  \\
\hline
\end{array}
\)

\end{center}

\medskip

\(\mikH=w-b = 2h\)  
\[ \mykappa{2h-2} = - (2h-1) \alpha\;,\quad 
 \mykappa{2h-3} = -3 (2h-1) \alpha\]
\[ \begin{vmatrix} \mykappa{2h-2} & H (H-1) \\
2\alpha^2 & -\alpha (H-1) \end{vmatrix}
= \begin{vmatrix} -2(2h-1) \alpha & 2h (2h-1) \\
2\alpha^2 & -\alpha (2h-1)  \end{vmatrix} = - \alpha^{2} (2h-1)(2h+1)\]

\[ \begin{vmatrix} \mykappa{2h-3} & (H-1) (H-2) \\
6\alpha^2 & -3 \alpha (H-1) \end{vmatrix}
= \begin{vmatrix} -3(2h-1) \alpha & (2h-1) (2h-2) \\
6\alpha^2 & - 3 \alpha (2h -1) \end{vmatrix} = 3 \alpha ^{2} (2h-1)(2h+1)\]

Thus, 
\begin{align*}
& \{ \myTF{1}{2a,b,H-2a} \mid 0\leqq 2a \leqq H\} \quad \text{are linearly
independent, and }\\
& \{ \myTF{1}{2a+1,b,H-1-2a} \mid 1\leqq 2a+1 \leqq H\} \quad \text{are linearly
independent, and so }\\
& \{ \myTF{1}{a,b,2h-a} \mid 0\leqq a \leqq 2h \} \quad \text{are linearly
independent. }
\\
& \{ \myWE{1}{a,b,2h-a} \mid 0\leqq a \leqq 2h \} \quad \text{are linearly
independent. }
\end{align*}

By the above discussions, 
we have a result about \( \dim \{ \myWE{1} {a,b,c} \mid a+c = w - b \}\),  which is similar to  Lemma \ref{lemma:slant:mu}.   
We denote  \( \Hori{w} {b} = \{ \myWE{1} {a,b,c} \mid a+c = w - b \}\).  
\begin{lemma}
\label{lemma:horizontal:L} 
If \(w-b=\text{ even}\),   
\(\dim \Hori{w}{b} = {}^{\#}\mathfrak{S}^{[w]}_{b} = w-b+1\),  (i.e., full rank).  

If \(w-b=\text{ odd}\),  then 
\(\dim \Hori{w}{b} = 
{}^{\#}\Hori{w}{b} - 1 = w - b \) and  
\( \{ \myWE{1}{a,b,c} \in \Hori{w}{b} 
\mid  c \ne 1 \}\) are  a base, by which 
\( \myWE{1}{w-b-1,b,1}\) are expressed.    
\end{lemma}
As a direct result, we have 
\begin{kmCor} \label{lemma:w+2:EW1}
\(\ds \dim \{ \myWE{1} {a,b,c} \mid a+b+c = w \} = 
\binom{w+2}{2} - \frac{w+1-\epsilon}{2}\), where
\( \epsilon = 1\) if  \(w\) is even, \(\epsilon = 0\)  otherwise.  
\end{kmCor}
\textbf{Proof:} From corollary above, 
\begin{align*}
\dim \{ \myWE{1}{a,b,c} \mid a+b+c=w \} &= \sum_{b=0}^{w} \Hori{w}{b} = 
\sum_{b=0}^{w} (w-b) + \# \{ w-b =\text{ even } \mid 0\leqq b \leqq  w\} \\
&= w(w+1) - \frac{w(w+1)}{2} +  \# \{ 2b  \mid 0\leqq b \leqq w\} \\
&=  \frac{w(w+1)}{2} +  \frac{w+1+\epsilon}{2} 
=   \frac{(w+2)(w+1)}{2} -  \frac{w+1-\epsilon}{2} \; .
\end{align*}
The parity of \(w\) still has some effect.  \( \epsilon = 2[w/2] -w +1\)
holds  if we use the Gauss symbol. 
\kmqed 

So far, we got 
\(\binom{w}{2} - \epsilon, \binom{w+2}{2} \) by \(\MC{}{P}\) or \(
\LC{3}{A}\), and 
\( \binom{w+2}{2} - \frac{w+1-\epsilon}{2}\) by 
Lemma \ref{lemma:w+2:EW1} generators.

The last thing to do is  to find which \(\myWE{2}{a,b,c} \) contributes
to enlarging the linearly independent system.  By bare hands, or by an experiment of symbol calculus, 
we know that 
 \(\{\myWE{2}{2a,w-2a,0} \mid  2a \leqq w\}\)
contribute as linear independent members, where   
\begin{align*} \myWE{2}{2a,w-2a,0}
=& 
2 \alpha \beta  \LC{1}{2 a-1,w-2 a-1,2}
-\alpha (2 a-w-1)  \LC{1}{2 a-1,w-2 a+1 ,0}
-2 \beta (2 a+1)  \LC{1}{2 a+1,w-2 a-1,0}
\\& 
-2 \beta^2  \LC{2}{2 a,w-2 a-2,2}
-\beta (2 a-w+1)  \LC{2}{2 a,w-2 a,0} \; .  \end{align*}
\kmcomment{
 \(\{\myWE{2}{2a, w-2a,0} \mid  2a \leqq w\}\) 
are candidate for linear independent systems. 
}
The cardinality is \(\frac{ w+1+\epsilon }{2}\).  This means the rank
\(\text{Rk}_{[w+2]}\) should be \begin{align*}
\text{Rk}_{[w+2]} & \geqq 
(\binom{w}{2} - \epsilon)  +  \binom{w+2}{2} 
+ (\binom{w+2}{2} - \frac{w+1-\epsilon}{2}) 
+ \frac{ w+1+\epsilon }{2}
\\& = \binom{w}{2} + 2 \binom{w+2}{2}
\end{align*}
Since the Betti number is given by  
\(\text{Bet}_{[w+2]} = \dim \wtedC{w+2}{w} -  \text{Rk}_{[w+2]} - \text{Rk}_{[w+3]}  \) and
it is non-negative,  
\(
 \text{Rk}_{[w+2]} \leqq 
\dim \wtedC{w+2}{w} -   \text{Rk}_{[w+3]}  = \binom{w}{2} + 2 \binom{w+2}{2}\)
holds in our case. Thus, we conclude   
the rank is
\(  \text{Rk}_{[w+2]} = \binom{w}{2} + 2 \binom{w+2}{2}\) and so the 
kernel dimension of \(\wtedC{w+3}{w}\)  is \(\binom{w+2}{2}\). 


%

\paragraph{The space of cycles in \(\wtedC{w+1}{w}\):} 
When \(w=0\), then \(\wtedC{w+1}{w} = \wtedC{1}{0} = \frakg\) and the
boundary operator is trivial, and so the kernel dimension is \(\dim
\frakg\). We then 
study the space of cycles for \(w > 0\).

Take a general chain \(L+M\) where \( L = \sum \LC{i}{a,b,c}  \zb{i}
\mywedge \bU{a,b,c}{w}\) and \( M = \sum \MC{i} {p,q,r} W[\eps{i}=0 ]
\mywedge \bU{p,q,r}{w-2} \) with unknown scalars \(\LC{i}{a,b,c}, \MC{i}{p,q,
r}\).   
\begin{align*}
\pdel(L+M) =& \sum \LC{i}{A} (
- \zb{i}\mywedge \pdel\bU{A}{w} + \SbtES{\zb{i}}{\bU{A}{w}}) \\&
+ \sum \MC{j}{P} (( \pdel W[\eps{j}=0] ) \mywedge \bU{P}{w-2} 
+ \SbtES{
 W[\eps{j}=0] 
}{\bU{P}{w-2}} )
\\
=& - \sum \zb{i} \mywedge \zb{4} \mywedge U^{p,q,r} (   
 \frac{2}{\alpha} \tbinom{p+2}{2} \LC{i}{p+2,q,r}
+ \frac{2}{\beta} \tbinom{q+2}{2} \LC{i}{p,q+2,r}
+ 2 \tbinom{r+2}{2} \LC{i}{p,q,r+2}
)
+\sum \LC{i}{A} \SbtES{ \zb{i} } {U^{A}}
\\& 
+ \zb{1} \mywedge \zb{4}\mywedge \sum \alpha  \MC{1}{P} U^{P}
+ \zb{2} \mywedge \zb{4}\mywedge\sum (-\beta ) \MC{2}{P} U^{P}
+ \zb{3} \mywedge \zb{4}\mywedge\sum  \MC{3}{P} U^{P}
\\& 
+ \sum \MC{1}{P} ( 
  \zb{3} \mywedge \zb{4} \mywedge \SbtES{ \zb{2} } {U^{P}}
- \zb{2} \mywedge \zb{4} \mywedge \SbtES{ \zb{3} } {U^{P}}
)
\\& 
+ \sum \MC{2}{P} ( 
  \zb{3} \mywedge \zb{4} \mywedge \SbtES{ \zb{1} } {U^{P}}
- \zb{1} \mywedge \zb{4} \mywedge \SbtES{ \zb{3} } {U^{P}}
)
\\& 
+ \sum \MC{3}{P} ( 
  \zb{2} \mywedge \zb{4} \mywedge \SbtES{ \zb{1} } {U^{P}}
- \zb{1} \mywedge \zb{4} \mywedge \SbtES{ \zb{2} } {U^{P}}
) 
\\ 
=&  
+\sum \LC{i}{A} \SbtES{ \zb{i} } {U^{A}}
\\& 
+ \zb{1} \mywedge \zb{4} \mywedge( \sum  \alpha \MC{1}{P} U^{P} 
- \sum \LC{1}{A} (\pdel U^{A}) /\zb{4}  
- \sum \MC{2}{P} \SbtES{ \zb{3} } {U^{P}}
- \sum \MC{3}{P} \SbtES{ \zb{2} } {U^{P}})
\\& 
+ \zb{2} \mywedge \zb{4} \mywedge ( \sum -\beta \MC{2}{P} U^{P} 
- \sum \LC{2}{A} (\pdel U^{A}) /\zb{4}  
- \sum \MC{1}{P} \SbtES{ \zb{3} } {U^{P}}
+ \sum \MC{3}{P} \SbtES{ \zb{1} } {U^{P}})
\\&
+ \zb{3} \mywedge \zb{4} \mywedge ( \sum  \MC{3}{P} U^{P} 
- \sum \LC{3}{A} (\pdel U^{A}) /\zb{4}  
+ \sum \MC{1}{P} \SbtES{ \zb{2} } {U^{P}}
+ \sum \MC{2}{P} \SbtES{ \zb{1} } {U^{P}})
\end{align*} 

\kmcomment{
\\ \sum_{A} \LC{i}{A} \SbtES{\zb{1}}{U^{A}} 
&= \sum_{A}

((b+1)\LC{i} {a,b+1,c-1}-\beta (c+1)\LC{i}{a,b-1,c+1}) 

U^{A}\;,
\\ \sum_{A} \LC{i}{A} \SbtES{\zb{2}}{U^{A}} 

&= \sum_{A}
(- (a+1)\LC{i} {a+1,b,c-1}+\alpha (c+1)\LC{i}{a-1,b,c+1}) 

U^{A}\;,

\\ \sum_{A} \LC{i}{A} \SbtES{\zb{3}}{U^{A}} 
&= \sum_{A}
( \beta (a+1)\LC{i} {a+1,b-1,c}-\alpha (b+1)\LC{i}{a-1,b+1,c}) 
U^{A}\;,
}

Thus, the linear equation \( \pdel (L+M) =0\) is defined by 
\begin{subequations}
\begin{align}
G_{0} =&  \label{G3:D3:type1:w+1:0}
(b+1)\LC{1} {a,b+1,c-1}-\beta (c+1)\LC{1}{a,b-1,c+1} 
- (a+1)\LC{2} {a+1,b,c-1}+\alpha (c+1)\LC{2}{a-1,b,c+1} 
\\&\quad 
+  \beta (a+1)\LC{3} {a+1,b-1,c}-\alpha (b+1)\LC{3}{a-1,b+1,c}
\notag
\\
G_{1} =& \label{G3:D3:type1:w+1:1}
\alpha \MC{1}{p,q,r} 
-\frac{1}{\alpha} 2\tbinom{p+2}{2} \LC{1}{p+2,q,r}
-\frac{1}{\beta} 2\tbinom{q+2}{2} \LC{1}{p,q+2,r}
- 2\tbinom{r+2}{2} \LC{1}{p,q,r+2}
\\& 
- \beta (p+1)\MC{2} {p+1,q-1,r}+\alpha (q+1)\MC{2}{p-1,q+1,r} 
\notag
\\& 
+  (p+1)\MC{3}{p+1,q,r-1}-\alpha (r+1)\MC{3}{p-1,q,r+1} 
\notag 
\\
G_{2} =& \label{G3:D3:type1:w+1:2}
- \beta \MC{2}{p,q,r}  
-\frac{1}{\alpha} 2\tbinom{p+2}{2} \LC{2}{p+2,q,r}
-\frac{1}{\beta} 2\tbinom{q+2}{2} \LC{2}{p,q+2,r}
- 2\tbinom{r+2}{2} \LC{2}{p,q,r+2}
\\& 
-  \beta (p+1)\MC{1} {p+1,q-1,r}+\alpha (q+1)\MC{1}{p-1,q+1,r} 
\notag
\\&
+ (q+1)\MC{3} {p,q+1,r-1}-\beta (r+1)\MC{3}{p,q-1,r+1} 
\notag
\\
G_{3} = & \label{G3:D3:type1:w+1:3}
\MC{3}{p,q,r} 
-\frac{1}{\alpha} 2\tbinom{p+2}{2} \LC{3}{p+2,q,r}
-\frac{1}{\beta} 2\tbinom{q+2}{2} \LC{3}{p,q+2,r}
- 2\tbinom{r+2}{2} \LC{3}{p,q,r+2}
\\&
- (p+1)\MC{1} {p+1,q,r-1}+\alpha (r+1)\MC{1}{p-1,q,r+1} 
\notag
\\& 
+ (q+1)\MC{2} {p,q+1,r-1}-\beta (r+1)\MC{2}{p,q-1,r+1} 
\notag
\end{align} 
\end{subequations}
If we focus on  \(\MC{3}{p,q,r}\) in \eqref{G3:D3:type1:w+1:3} and
substitute them to  \eqref{G3:D3:type1:w+1:1} and \eqref{G3:D3:type1:w+1:2},  
then we have next equations without \( \MC{1}{P}\) and \( \MC{3}{P}\): 
\begin{align}
H_{1} (p,q,r) & =  
 -(r+2) (r+1)  \LC{1}{p,q,r+2}  
 -(p+2) (p+1) \frac{1}{\alpha } \LC{1}{p+2,q,r}  
 -(q+2)  (q+1)\frac{1}{\beta}   \LC{1}{p,q+2,r}
\\& 
+ (p+3) (p+2) (p+1) \frac{1}{\alpha}   \LC{3}{p+3,q,r-1}  
+ (q+2)  (q+1) (p+1)\frac{1}{\beta}   \LC{3}{p+1,q+2,r-1}  
\notag
\\& 
 -(r+1) (p+1) (p-r)   \LC{3}{p+1,q,r+1}  
 - \alpha (r+3) (r+2) (r+1)   \LC{3}{p-1,q,r+3}  
 \notag
 \\& 
 -\alpha (q+2) (q+1) (r+1)\frac{1}{\beta}   \LC{3}{p-1, q+2,r+1}
+ \alpha (r+2) (q+1)   \MC{2}{p-1,q+1,r}  
\notag
\\& 
+ \beta (r-1) (p+1)   \MC{2}{p +1,q-1,r}  
 -(q+1) (p+1)   \MC{2}{p+1,q+1,r-2}  
 \notag\\
 & 
 -\beta \alpha (r+2) (r+1)   \MC{2}{p-1,q-1,r+2} \;, 
 \notag
\\
H_{2} (p,q,r) & =  
 -(r+2) (r+1)    \LC{2}{p,q,r+2} 
 -(p+2) (p+1)\frac{1}{\alpha }   \LC{2}{p+2,q,r} 
 -(q +2) (q+1)\frac{1}{\beta}    \LC{2}{p,q+2,r}
 \\& 
+ (p+2) (p+1) (q+1)\frac{1}{\alpha}   \LC{3}{p+2,q+1,r-1}  
+ (q+ 3) (q+2) (q+1)\frac{1}{\beta}    \LC{3}{p,q+3,r-1}  
\notag
\\& 
 -(r+1) (q+1) (q-r)    \LC{3}{p,q+1,r+1}  
 -\beta  (r+3) (r+2) (r+1)    \LC{3}{p,q-1,r+3}  
 \notag
 \\& 
 -\beta (p+2) (p+1) (r+1)\frac{1}{\alpha}    \LC{3}{p+2,q-1,r +1}
+\beta (2 q r+q+r-1)    \MC{2}{p,q,r}  
\notag
\\& 
 -(q+2) (q+1)    \MC{2}{p,q+2,r-2} 
 -\beta^2 (r+2) (r+1)    \MC{2}{p,q-2,r+2} \;. 
 \notag
\end{align}
Now we look for the dimension of the space generated by \( G_{0}\), \( H_{1}
\) and \( H_{0}\), and we get \( \binom{w+2}{2} + \binom{w}{2}\).  
Thus, the rank is \( 2 \binom{w}{2} + \binom{w+2}{2}\) and   
the kernel dimension is \(  \binom{w}{2} + 2 \binom{w+2}{2}\).    

\kmcomment{

Oo := (p,q,r) -> 
- (
(p+1)* (q+1)* N[p+1,q+1,r] - (r+2)*(r+1)*   N[p,q,r+2]
 # + O[p,q,r] 
   -(p+1)* P[p+1,q,r-1]-2*(r+1)*P[p,q-1,r+1]
+   (q+1)* Q[p,q+1,r-1] +2*(r+1)*Q[p-1,q,r+1] )  ;

(* 
hh1 := {[(b+1)*(a+1), L[a+1,b+1,c]], [-(c+2)*(c+1), L[a,b,c+2]]}, {}, {[-(a+2)
*(a+1)*(b+1), N[a+2,b+1,c-1]], [-(c+1)*(-c+2*b)*(a+1), N[a+1,b,c+1]], [2*(c+3)
*(c+2)*(c+1), N[a,b-1,c+3]]}, {}, {};

hh2 := {}, {[(b+1)*(a+1), M[a+1,b+1,c]], [-(c+2)*(c+1), M[a,b,c+2]]}, {[-(b+2)
*(b+1)*(a+1), N[a+1,b+2,c-1]], [-(c+1)*(b+1)*(2*a-c), N[a,b+1,c+1]], [2*(c+3)*
(c+2)*(c+1), N[a-1,b,c+3]]}, {}, {};
*) 

H1 := (a,b,c) ->  (b+1)*(a+1)* L[a+1,b+1,c] -(c+2)*(c+1) * L[a,b,c+2]
-(a+2) *(a+1)*(b+1)*N[a+2,b+1,c-1]-(c+1)*(-c+2*b)*(a+1)* N[a+1,b,c+1]
+2*(c+3) *(c+2)*(c+1)*N[a,b-1,c+3]:

H2 := (a,b,c) ->  
(b+1)*(a+1)* M[a+1,b+1,c]-(c+2)*(c+1)* M[a,b,c+2]
-(b+2) *(b+1)*(a+1)* N[a+1,b+2,c-1]-(c+1)*(b+1)*(2*a-c)* N[a,b+1,c+1] +
2*(c+3)* (c+2)*(c+1)* N[a-1,b,c+3] :

HH1 := (a,b,c) ->  (b+1)*(a+1)* L[a+1,b+1,c] -(c+2)*(c+1) * L[a,b,c+2]
-(a+2) *(a+1)*(b+1)*NN(a+2,b+1,c-1)-(c+1)*(-c+2*b)*(a+1)* NN(a+1,b,c+1)
+2*(c+3) *(c+2)*(c+1)*NN(a,b-1,c+3):

HH2 := (a,b,c) ->  
(b+1)*(a+1)* M[a+1,b+1,c]-(c+2)*(c+1)* M[a,b,c+2]
-(b+2) *(b+1)*(a+1)* NN(a+1,b+2,c-1)-(c+1)*(b+1)*(2*a-c)* NN(a,b+1,c+1) +
2*(c+3)* (c+2)*(c+1)* NN(a-1,b,c+3) : 

}

\paragraph{The space of cycles in \(\wtedC{w}{w}\):} 
When \(w=0\), then \(\wtedC{w}{w} = \Lambda ^{0}\frakg\) and the
boundary operator is trivial, and so the kernel dimension is 1.  
We study the space  of cycles in \(\wtedC{w}{w}\)  for \(w > 0\).  
Take a general chain \(L+M\) where \( L = \sum \LC{ }{A}
\bU{A}{w}\) and \( M = \sum \MC{i} {P}  \zb{i} \mywedge \zb{4}
\mywedge \bU{P}{w-2} \) with unknown scalars \(\LC{ }{A},
\MC{i}{P}\).   
Since 
\begin{align*}
\pdel(L+M) =&  \sum \LC{}{A} \pdel U^{A} 
+ 0 + 0 
+ \sum \MC{j}{P} \zb{4} \mywedge \SbtES{\zb{j}}{U^{P}}
\end{align*}
the linear equations \( \pdel(L+M)=0\) is defined by 
\begin{align}
& \label{HHK:T1:0:a}
\frac{(p+2)(p+1)}{\alpha} \LC{}{p+2,q,r}
+ \frac{(q+2)(q+1)}{\beta} \LC{}{p,q+2,r}
+ (r+2)(r+1) \LC{}{p,q,r+2}
\\& + (q+1) \MC{1} {p,q+1,r-1} - \beta (r+1) \MC{1}{p,q-1,r+1}
- (p+1) \MC{2}{p+1,q,r-1} + \alpha (r+1) \MC{2}{p-1,q,r+1}
\notag
\\& 
+ \beta (p+1) \MC{3}{p+1,q-1,r} - \alpha (q+1) \MC{3}{p-1,q+1,r}
\notag
\end{align}
\kmcomment{ \textcolor{red}{Revised} on July 21, 2020 by focusing the rank,
the number of generators of linear equations: }
Thus, the rank is \( \binom{w}{2}\), and the kernel dimension is \(
2\binom{w}{2} + \binom{w+2}{2}\).  

\kmcomment{
\paragraph{The space of cycles in \(\wtedC{w-1}{w}\):} 
Generators of this chain space are \( \zb{4}\mywedge \bU{A}{w-2}\). 

Since 
\(\pdel (\zb{4}\mywedge \bU{A}{w-2}) = 
- \zb{4}\mywedge \pdel(\bU{A}{w-2}) + \SbtES{\zb{4}}{ \bU{A}{w-2}} 
= - \zb{4} \mywedge \zb{4} \wedge  (1-\alpha ) b c \bU{a,b-1,c-1}{w-2} 
+ 0 = 0 \), the boundary operator is trivial and the kernel dimension is 
\( \binom{w}{2}\). 
}

\renewcommand{\arraystretch}{1.3}

\paragraph{Final table of chain complex }
\begin{thm} \label{thm:G3D3-T1} \     
\\
\begin{center}
\(
\begin{array}{c|*{5}{c}}
\text{weight}= w > 0 & w-1 & w & w+1 & w+2 & w+3\\\hline
\text{SpaceDim} & 
\binom{w}{2} & 
3 \binom{w}{2} + \binom{w+2}{2} & 
3\binom{w}{2} + 3\binom{w+2}{2} & 
\binom{w}{2} + 3\binom{w+2}{2} & 
\binom{w+2}{2} 
\\
\ker\dim & \binom{w}{2} & 2\binom{w}{2}  +\binom{w+2}{2}  & 
\binom{w}{2} + 2 \binom{w+2}{2}   & \binom{w+2}{2}  & 0   \\\hline
\text{Betti} & 0 & 0 &  0  & 0   &  0 
\end{array}
\)
\end{center}
\end{thm}

\section{Further questions}
So far, we discussed super homology groups associated with Lie algebras of
dimension 2 or 3.  We observed that if \( \frakg = \Sbt{\frakg}{\frakg}\)
then the Betti numbers are all zero, and it may be interesting question if
this is true when the dimension is larger than 3.  

As an example of 4-dimensional Lie algebra, we take \(\mathfrak{gl}(2,\mR)\) 
of \(2 \times 2\) matrices, which is an extension of 
\(\mathfrak{sl}(2,\mR)\). 
By experiments for lower weights, we have the following information about
Betti numbers (we do not note dimension or rank for each chain space here).   
\begin{center}
\setlength{\arraycolsep}{5pt}
\renewcommand{\arraystretch}{1.0}
\( 
\begin{array} { c | *{10}{r|} }
m\text{-th chain}  & 0 & 1 & 2 & 3 & 4 & 5 & 6 & 7 & 8 & 9 \\ \hline
w=0 & 1 & 1 & 0 & 1 & 1 \\ \hline
w=1 &   & 0 & 0 & 0 & 0 & 0 \\ \hline
w=2 &   & 0 & 2 & 2 & 0 & 2 & 2 \\ \hline
w=3 &   & 0 & 1 & 1 & 0 & 1 & 1 & 0 \\ \hline
w=4 &   &   & 0 & 0 & 2 & 2 & 0 & 2 & 2 \\ \hline
w=5 &   &   & 0 & 1 & 2 & 1 & 1 & 2 & 1 & 0 \\ \hline
\end{array}
\)
\end{center}

Looking at this table, we expect some ``rule''' and expect the rule comes
from that the Lie algebra is an extension of \( \mathfrak{sl}(2,\mR)\). 


\nocite{Mik:Miz:super2}
\nocite{Mik:Miz:super3}

\bibliographystyle{plain}
\bibliography{km_refs}
\end{document}